\documentclass[11pt]{amsart}
\usepackage{fullpage}
\def\showanswers{0}

\newcommand{\hide}[1]{
\ifnum\showanswers=1
        {\color{red} #1} %\vspace{\baselineskip}
        \fi
\ifnum\showanswers=0
        \fi
}

\usepackage{hyperref}
\hypersetup{
    colorlinks=false,
    linkcolor=,
    filecolor=,      
    urlcolor=,
    pdftitle={Overleaf Example},
    pdfpagemode=FullScreen,
    }

\usepackage{comment} 

\usepackage{enumitem}

\usepackage{mathrsfs}
\usepackage{amsthm}
\usepackage{amsmath}
\usepackage{amssymb}

\theoremstyle{plain}
\newtheorem{lemma}{Lemma}[section]
\newtheorem{sublemma}{Lemma}[lemma]
\newtheorem{proposition}[lemma]{Proposition}
\newtheorem{theorem}[lemma]{Theorem}
\newtheorem{corollary}[lemma]{Corollary}
\newtheorem{maintheorem}{Theorem}

\newtheorem*{keylemma*}{Key Lemma}

\theoremstyle{definition}
\newtheorem{question}[lemma]{Question}
\newtheorem*{openquestion*}{Open Question}
\newtheorem{definition}[lemma]{Definition}
\newtheorem{example}[lemma]{Example}

\newtheorem{remark}[lemma]{Remark}
\newtheorem*{remark*}{Remark}
\newtheorem{notation}[lemma]{Notation}

\newtheorem{setting}[lemma]{Setting}

\newcommand{\piv}{\pi_{\mathcal{V}}}
\newcommand{\piz}{\pi_{\mathcal{Z}}}

\usepackage[T1]{fontenc}
\usepackage{mathtools}

\newcommand{\Addresses}{{% additional braces for segregating \footnotesize
  \bigskip
  \hide{\footnotesize}

  \textsc{Department of Mathematics, University of Oklahoma, 601 Elm Ave, Norman, OK 73019-3103, USA}\par\nopagebreak
  \textit{Email address}: \texttt{tomoya.tatsuno@ou.edu}

}}

\title{Abstract Group Automorphisms of Two-Step Nilpotent Lie Groups and Partial Automatic Continuity}
\author{Tomoya Tatsuno}
\date{}

\begin{document}

\begin{abstract}
    Any abstract (not necessarily continuous) group automorphism of a simple, compact Lie group must be continuous due to Cartan (1930) and van der Waerden (1933). The purpose of this paper is to study a similar question in nilpotent Lie groups. For uncountably many simply connected 2-step nilpotent Lie groups, we show any abstract group automorphism is continuous ``up to discontinuity due to the center and ring automorphisms of $\mathbb{R}$-algebras.'' Such groups include (1) a generic simply connected 2-step nilpotent Lie group $N$ with $\dim [N,N]\geq \frac{\dim N-1}{2}$, (2) all the 12 simply connected 2-step nilpotent Lie groups of dimensions 6 or less, and (3) Iwasawa N-groups of simple Lie groups of rank 1, i.e., ``Heisenberg groups.'' All of the three cases are derived from one key result that gives a sufficient condition for the automorphism group to be of the type described above. Previously, only countably many simply connected nilpotent Lie groups were known to satisfy the same property. We obtain similar results for complex Lie groups as well. We also show if any Lie group $G$ and the group $N$ of (1) or (3) are isomorphic as abstract groups, then they are isomorphic as Lie groups. 
\end{abstract}

\maketitle
\Addresses

\subsubsection*{Keywords}
Nilpotent Lie Groups; Generic Lie Groups; Heisenberg Groups; Iwasawa Decomposition          

\subsubsection*{Mathematics Subject Classification} 22E25, 22E46, 20F28, 53C30 

\subsubsection*{Declarations of interest:} none 

\section{Introduction} 
\subsection{Abstract Group Automorphisms} \label{section_abstract_group_automorphisms}
The purpose of this paper is to obtain a complete description of abstract (not necessarily continuous) group automorphisms of a large class of 2-step nilpotent Lie groups. In particular, we present rigidity results.

Our motivation comes from a classical rigidity result called \textit{automatic continuity}: \textit{any abstract (not necessarily continuous) group isomorphism between two connected, compact simple Lie groups must be continuous} (independently proved by Cartan \cite{Car1930} and van der Waerden \cite{vanderwaerden1933}). In particular, any abstract group automorphism of such a Lie group is a Lie group automorphism. In 1941, Freudenthal showed that automatic continuity holds for certain non-compact simple Lie groups, including $SL_n(\mathbb{R})$ (\cite{Freudenthal_1941}). These results in the simple setting have been greatly generalized by many authors; see, e.g., \cite{Braun_Hofmann_Kramer}, and references therein.

However, automatic continuity holds only for very special Lie groups. Even in simple Lie groups, discontinuous group automorphisms exist. Take the group $SL_2(\mathbb{C})$ and a discontinuous field automorphism $\varphi$ of the complex numbers $\mathbb{C}$. Define $\hat{\varphi}:SL_2(\mathbb{C})\to SL_2(\mathbb{C})$ by
\begin{align*}
    \hat{\varphi}\left(\begin{array}{cc}
        a & b \\
        c & d
    \end{array}\right)=\left(\begin{array}{cc}
        \varphi(a) & \varphi(b) \\
        \varphi(c) & \varphi(d)
    \end{array}\right).
\end{align*}
In the literature, the group automorphism $\hat{\varphi}$ is called a \textit{field automorphism} of $SL_2(\mathbb{C})$. This is a discontinuous group automorphism of $SL_2(\mathbb{C})$, first observed by von Neumann. For $SL_2(\mathbb{C})$, a composition of any Lie group automorphism and $\hat{\varphi}$ also gives rise to a discontinuous group automorphism.

It turns out that this is the only way to get discontinuous group automorphisms of $SL_2(\mathbb{C})$. More precisely, any group automorphism $F$ of $SL_2(\mathbb{C})$ has the form $F = \overline{F} \circ \hat{\varphi}$, where $\overline{F}$ is a Lie group automorphism and $\hat{\varphi}$ is a field automorphism  of $SL_2(\mathbb{C})$ induced by an automorphism $\varphi$ of the field $\mathbb{C}$ (such automorphisms are called ``standard''). These kinds of results were proved in much greater depth by Borel and Tits in 1974 for abstract group homomorphisms between certain simple algebraic groups over infinite fields; see \cite{Borel_and_Tits}.

Little is known about abstract group automorphisms of Lie groups in general. In this article, we shed light on the rigidity of abstract automorphisms of certain nilpotent Lie groups. It was Kallman and McLinden \cite{Kallman_McLinden} who first pointed out that the 3-dimensional Heisenberg group admits a discontinuous group automorphism $\mu:\textrm{Heis}_3(\mathbb{R})\to \textrm{Heis}_3(\mathbb{R})$ of the following form:

\begin{align*}
    \mu\left(\begin{array}{cccc}
            1 & x & z \\
            0 & 1 & y\\
            0 & 0 & 1 
        \end{array}\right) = \left(\begin{array}{cccc}
            1 & x & z+\lambda(x,y) \\
            0 & 1 & y\\
            0 & 0 & 1 
        \end{array}\right),
\end{align*}
where $\lambda:\mathbb{R}^2\to \mathbb{R}$ is a discontinuous $\mathbb{Q}$-linear map (see \cite[Corollary 5.2.2]{Kuczma}). Moreover, they showed that any abstract group isomorphism from any Polish group to $\textrm{Heis}_3(\mathbb{R})$ possesses partial automatic continuity (see \cite[Proposition 5.1]{Kallman_McLinden}). 

We point out that the automorphism $\mu$ above is precisely a \textit{central automorphism}, that is,  a centralizer of the group of inner automorphisms. Central automorphisms form a normal subgroup of the group of abstract group automorphisms. Any central automorphism of $\textrm{Heis}_3(\mathbb{R})$ has the form $\mu$ above (with $\lambda$ possibly discontinous). Moreover, we show any simply connected nilpotent Lie group admits a discontinuous central automorphism, and it comes from $\mathbb{Q}$-linear maps (see Section \ref{section_central_automorphisms}).

Another interesting class of discontinuous group automorphisms arises from discontinuous ring automorphisms of $\mathbb{R}$-algebras, first observed by Tits \cite{Tits_Lie}. Let $\mathbb{R}[\epsilon]=\mathbb{R}[x]/\langle X^2\rangle$ be the ring of dual numbers, where $\epsilon$ denotes $X$ in the quotient. There exists a discontinuous ring automorphism $\varphi$ of $\mathbb{R}[\epsilon]$ (see \cite[2.4]{Tits_Lie}), and it induces a discontinuous group automorphism $\hat{\varphi}$ of $\textrm{Heis}_3(\mathbb{R}[\epsilon])$ by applying $\varphi$ to each entry, as in the case of $SL_2(\mathbb{C})$ above. We call the group automorphism $\hat{\varphi}$ a \textit{ring automorphism} of $\textrm{Heis}_3(\mathbb{R}[\epsilon])$. Observe, this generalizes discontinuous group automorphisms induced by field automorphisms of $\mathbb{C}$, as $\mathbb{C}$ is also a $\mathbb{R}$-algebra.

Any composition of a discontinuous central automorphism, a Lie group automorphism, and a ring automorphism gives rise to a discontinuous group automorphism of a simply connected nilpotent Lie group. Naturally, it raises the question of whether or not this is the only possible way for the discontinuity to occur.

\begin{question}\label{driving_question}
    Let $N$ be a simply connected nilpotent Lie group. Let $F$ be an abstract (not necessarily continuous) group automorphism of $N$. For which $N$, does $F$ always have the form
    \begin{align*}
        F = \mu \circ \overline{F}\circ \hat{\varphi},
    \end{align*}
    where $\mu$ is a central automorphism, $\overline{F}$ is a Lie group automorphism, and $\hat{\varphi}$ is a ring automorphism?
\end{question}

Following the terminology of Kallmann--McLinden \cite[Section 5]{Kallman_McLinden}, we give a definition.
\begin{definition}\label{def_of_partial_automatic_continuity}
    Let $N$ be a simply connected nilpotent Lie group. If any group automorphism $F$ has the form as in Question \ref{driving_question}, we say that $N$ satisfies \textit{partial automatic continuity}.
\end{definition}
\begin{remark}
    It follows from \cite{Kallman_McLinden} that $\textrm{Heis}_3(\mathbb{R})$ satisfies partial automatic continuity.
\end{remark}
\begin{remark}
    In \cite{Kallman_McLinden}, the term partial automatic continuity was used for abstract group \textit{isomorphisms} from a Polish group to $\textrm{Heis}_3(\mathbb{R})$.
\end{remark}

One can deduce from classical results \cite{Gibbs}, \cite{Levchuk}, \cite{Khor}, \cite{Im_symplectic}, \cite{Im_typeD}, \cite{Myasnikov_Sohrabi} that countably many distinguished simply connected nilpotent Lie groups with a rational structure satisfy partial automatic continuity; see the ``History'' below. Our first result covers uncountably many generic groups without any rational structure, showing that partial automatic continuity holds beyond classic examples. 

\begin{maintheorem}\label{generic_condition_introduction}
A generic simply connected 2-step nilpotent Lie group $N$ with $\dim [N,N]\geq \frac{\dim N-1}{2}$ satisfies partial automatic continuity. For any abstract (not necessarily continuous) group automorphism $F$ of $N$, there exist a central automorphism $\mu$ and a Lie group automorphism $\overline{F}$ of $N$ such that 
    \begin{align*}
        F=\mu \circ \overline{F}.
    \end{align*}
\end{maintheorem}

\begin{remark}
    See Remark \ref{remark_reason_for_bounds} for the reason for the dimension bound. Generically, these groups do not have an additional underlying $\mathbb{R}$-algebra structure such as a complex structure, so one does not see a ring automorphism above.
\end{remark}

\begin{remark}
    The word generic refers to a Zariski open dense set of the space of the structure constants of the corresponding Lie algebras. See Theorem \ref{generic_condition_precise} in Section \ref{section_generic} for a more precise statement.
\end{remark}

\begin{remark}\label{why_uncountable}
    The moduli space of isomorphism classes has a positive dimension with a few exceptions (see \cite[Table 1]{Jablonski_moduli}), so one deduces that there are uncountably many non-isomorphic groups that satisfy partial automatic continuity. Such groups usually do not admit any rational structure. Even further, the set of algebras with a rational structure has measure zero (when the dimension is large enough).  
\end{remark}

Let $\dim[N,N]\geq \frac{\dim N-1}{2}$, $p=\dim[N,N]$, and $q=\dim N-p$. With a few exceptions, the group of abstract group automorphisms of a generic group $N$ is isomorphic to 
    \begin{align*}
        \{\lambda:\mathbb{R}^q\to \mathbb{R}^p\;|\;\textrm{$\lambda$ is $\mathbb{Q}$-linear}\} \rtimes K,
    \end{align*}
where $K$ is a one-dimensional Lie group. In a sense, this means a generic 2-step nilpotent Lie group $N$ has the smallest symmetry (as we always have $\dim K\geq 1$ in general). See the end of Section \ref{section_semidirect_product} for more details.

To our knowledge, Theorem \ref{generic_condition_introduction} is the first result concerning abstract group automorphisms of generic nilpotent Lie groups.  Whether or not this result holds for $\dim [N, N]  < \frac{\dim N-1}{2}$ remains open. 

\begin{question}
    Does any generic simply connected 2-step nilpotent Lie group satisfy partial automatic continuity?
\end{question}

We identify a convenient criterion that forces partial automatic continuity.

\begin{maintheorem} \label{maintheorem_criterion}
    Let $\mathcal{N}$ be a real 2-step nilpotent Lie algebra. Let $\mathcal{Z}$ be the center of $\mathcal{N}$. Let $N$ be a simply connected 2-step nilpotent Lie group with Lie algebra $\mathcal{N}$.

    If there exists a vector space complement $\mathcal{V}$ of $\mathcal{Z}$ in $\mathcal{N}$ and a basis $e_1,...,e_q$ of $\mathcal{V}$ such that $[e_1,. e_j]$, $j\in \{2, ..., q\}$, is linearly independent, then $N$ satisfies partial automatic continuity.
    
    For any abstract (not necessarily continuous) group automorphism $F$ of $N$, there exist a central automorphism $\mu$ and a Lie group automorphism $\overline{F}$ of $N$ such that 
    \begin{align*}
        F=\mu \circ \overline{F}.
    \end{align*}
\end{maintheorem}

See Corollary \ref{corollary_concatenation} for the most general case covered in this paper. Intuitively, if there are many Lie bracket relations, then one can hope for a more rigid structure. For example, a free 2-step nilpotent Lie algebra generated by 
$x_1, ..., x_m$ satisfies the assumption of Theorem \ref{maintheorem_criterion}.  

To connect Theorem \ref{maintheorem_criterion} with Theorem \ref{generic_condition_introduction}, observe that the linear independence condition in Theorem \ref{maintheorem_criterion} is a Zariski open condition. 

\begin{remark} \label{remark_reason_for_bounds}
    The assumption of Theorem \ref{maintheorem_criterion} implies $\dim [N,N] \geq \dim \mathcal{V}-1$. On the other hand, for a generic group, the center coincides with $[N,N]$, so $\dim \mathcal{V}=\dim N-\dim [N,N]$. This implies $\dim [N,N]\geq \frac{\dim N-1}{2}$. This is the reason for the dimension bound in Theorem \ref{generic_condition_introduction}.
\end{remark}

As an application of Theorem \ref{maintheorem_criterion} and a more general result (Corollary \ref{corollary_concatenation}), we obtain the following theorem.

\begin{maintheorem} \label{pac_up_to_dim_6_introduction}
    All the 12 simply connected 2-step nilpotent Lie groups of dimensions 6 or less satisfy partial automatic continuity.
\end{maintheorem}

See Section \ref{section_counterexamples} for more details. The classification is summarized in Table \ref{table_low_dimensions}.

\vskip\baselineskip
\noindent
\textbf{History.} Abstract group automorphisms of simply connected nilpotent Lie groups are implicitly studied in group theory. We explain how one can deduce that certain nilpotent Lie groups satisfy partial automatic continuity from classical results. We list all the known examples of simply connected nilpotent Lie groups that satisfy partial automatic continuity to the best of our knowledge.

Given a simple Lie group $G$ and an Iwasawa decomposition $G=KAN$, call the simply connected nilpotent Lie group $N$ the \textit{Iwasawa N-group} of $G$. 

In 1970, Gibbs studied abstract group automorphisms of maximal unipotent subgroups of Chevalley groups and twisted Chevalley groups over a field of characteristic $\neq 2,3$ \cite[Theorem 6.2, Theorem 7.1]{Gibbs}. From their work, one deduces that Iwasawa N-groups of complex simple Lie groups (viewed as a real Lie group) and split and quasi-split real simple Lie groups satisfy partial automatic continuity. Except for the 3-dimensional Heisenberg group, these groups are of 3-step or higher. 

In 1985, Khor showed a similar result for the unipotent radical of any parabolic subgroup of $GL(l, K)$ that is not the Borel subgroup for $l\geq 4$ and for a field $K$ of characteristic $\neq 2$  \cite{Khor}. Partial automatic continuity holds for these groups over $K=\mathbb{R}$ or $\mathbb{C}$. In particular, the 2-step nilpotent Lie groups covered in Khor's paper are of the block form \begin{align*}
 \left(\begin{array}{ccc}
            I & A_{12} & A_{13}\\
            0 & I & A_{23}\\
            0 & 0 & I
        \end{array}\right).
\end{align*}
The matrices $A_{12}, A_{13}$ and $A_{23}$ are block matrices of possibly different size. The diagonals are the identity square matrices.  Such a group contains a large abelian subgroup corresponding to $\langle A_{12}, A_{13}\rangle$ and $\langle A_{23}, A_{13}\rangle$. For this reason, except for the $(2n+1)$-dimensional Heisenberg groups with $n\geq 2$, all the 2-step nilpotent Lie groups covered by Khor are not \textit{nonsingular} (\cite[p.618]{Eberlein}). 

Im showed a similar result for the unipotent radical of certain parabolic subgroups of Chevalley groups of type C (\cite{Im_symplectic}, \cite{Im_dissertation}) and D (\cite{Im_typeD}). They assume that their flags have a length of 4 or more, so the nilpotent Lie groups treated there are of 5-step or more.

Other groups satisfying partial automatic continuity include $N$ whose Lie algebra $\mathcal{N}$ is a free nilpotent Lie algebra over any $\mathbb{R}$-algebra, or the algebra of strictly upper triangular matrices over any $\mathbb{R}$-algebra. See \cite{Myasnikov_Sohrabi} and \cite{Levchuk} respectively, for their more general work over rings.

As far as we know, these groups exhaust all the known examples of nilpotent Lie groups satisfying partial automatic continuity.

Among uncountably many non-isomorphic nilpotent Lie groups, one may consider Iwasawa N-groups as model groups. As Gibbs showed, the split case and quasi-split case satisfy partial automatic continuity. Therefore, it is natural to ask the following.

\begin{question}
    Do all Iwasawa N-groups satisfy partial automatic continuity?
\end{question}

The split real simple Lie groups have the highest possible real rank. We study the other extreme - the Iwasawa N-group $N$ of a real simple Lie group of rank 1, which is abelian or of 2-step. In this case, $N$ is one of the abelian Lie groups $\mathbb{R}^n$,  $(2n+1)$-dimensional Heisenberg groups, $(4n+3)$-dimensional quaternionic Heisenberg groups, and 15-dimensional octonionic Heisenberg group. The group $N$ plays an important role in the geometry of symmetric spaces of rank 1; see, e.g., \cite[pp. 617--618]{Eberlein}.

Partial automatic continuity is satisfied by the abelian Lie groups $\mathbb{R}^n$ (as any automorphism is central) and $(2n+1)$-dimensional Heisenberg groups (\cite[$n=1$]{Gibbs}, \cite[$n\geq 2$]{Khor}). However, until now, the abstract group automorphism group of the other two cases - quaternionic and octonionic - is not known. We now present our fourth main result.  

\begin{maintheorem} \label{temp_theorem1_intro}
    Let $N$ be the Iwasawa N-group of a real simple Lie group of rank 1, i.e., one of the abelian Lie group $\mathbb{R}^n$,  $(2n+1)$-dimensional Heisenberg groups, $(4n+3)$-dimensional quaternionic Heisenberg groups, and 15-dimensional octonionic Heisenberg group. Then, $N$ satisfies partial automatic continuity. For any abstract (not necessarily continuous) group automorphism $F$ of $N$, there exist a central automorphism $\mu$ and a Lie group automorphism $\overline{F}$ of $N$ such that 
    \begin{align*}
        F=\mu \circ \overline{F}.
    \end{align*}
\end{maintheorem}
\begin{remark} \label{remark_on_homomorphism}
    Theorem \ref{temp_theorem1_intro} does not extend to group homomorphisms; see Example \ref{discontinuous_homomorphism_preserving_v}. 
\end{remark}

The group of group automorphisms of $N$ has a semidirect product structure; see Section \ref{section_semidirect_product} and Remark \ref{remark_on_semidirect_product_structure_at_group_level}.

The base case $n=1$ is covered by Theorem \ref{maintheorem_criterion}, and the case $n>1$ follows from a more general result (Corollary \ref{corollary_concatenation}).

\subsection{Isomorphisms between Lie Groups}

We also consider the following question. We assume Lie groups are connected.
\begin{question}
    If Lie groups $G$ and $H$ are isomorphic as abstract groups, are they isomorphic as Lie groups?
\end{question}
If any abstract group isomorphism from $G$ to $H$ must be continuous, it follows that $G$ and $H$ are isomorphic as Lie groups.

On the other hand, $\mathbb{R}$ and $\mathbb{R}^2$ are isomorphic as abstract groups, but they are not isomorphic as Lie groups. Hence, $G=\textrm{Heis}_3(\mathbb{R})\times \mathbb{R}$ and $H=\textrm{Heis}_3(\mathbb{R})\times \mathbb{R}^2$ are isomorphic as abstract groups, but not as Lie groups. In this example, the issue comes from the abelian factor. Note, a simply connected 2-step nilpotent Lie group $N$ does not have an abelian factor if and only if $[N,N]=Z(N)$, the center. Moreover, any simply connected 2-step nilpotent Lie group is a direct product of a simply connected 2-step nilpotent Lie group $N$ with $[N,N]=Z(N)$ and $\mathbb{R}^k$ for some $k$. Therefore, it is reasonable to restrict our attention to the group $N$ with $[N,N]=Z(N)$.

\begin{maintheorem} \label{isomorphism_theorem_introduction}
    Let $N$ be a simply connected 2-step nilpotent Lie group in Theorem \ref{generic_condition_introduction} or Theorem \ref{temp_theorem1_intro}. Let $G$ be a Lie group. If $G$ and $N$ are isomorphic as abstract groups, then they are isomorphic as Lie groups. 
\end{maintheorem}
A generic group $N$ does not have an abelian factor, and Heisenberg groups do not have an abelian factor. A group $N$ in Theorem \ref{maintheorem_criterion} is allowed to have an abelian factor, and low-dimensional examples do have an abelian factor, so the groups $N$ in Theorem \ref{maintheorem_criterion} and Theorem \ref{pac_up_to_dim_6_introduction} are omitted here.

\begin{remark}
    More generally, this theorem holds for a concatenation of groups (without an abelian factor) covered in Theorem \ref{maintheorem_criterion} (see Definition \ref{definition_concatenation}). For the statement in full generality and the proof, see Theorem \ref{theorem_isomorphism_abstract_and_Lie_group}. Theorem \ref{isomorphism_theorem_introduction} does not extend to an affine algebraic group setting; see Example \ref{example_arising_from_field_of_rational_functions}.
\end{remark}

\subsection{Algebraic Groups} \label{section_introduction_complex_case}
Any simply connected nilpotent Lie groups are affine algebraic groups over $\mathbb{R}$, and some of our results hold true beyond the field $\mathbb{R}$. We discuss the case of simply connected complex 2-step nilpotent Lie groups, i.e., the case of affine algebraic groups over $\mathbb{C}$. See Remark \ref{remark_affine_algebraic_groups} for the case of a field $\mathbb{F}$ of characteristic not 2.

Given a complex Lie algebra $\mathcal{N}$ and a field automorphism $\varphi$ of $\mathbb{C}$, the $\varphi$-conjugate Lie algebra $\mathcal{N}^{\varphi}$ of $\mathcal{N}$ is defined as follows (\cite[Section 3]{Dere}). As vector spaces, $\mathcal{N}=\mathcal{N}^{\varphi}$, and if $X_1, ..., X_n$ is a $\mathbb{C}$-basis of $\mathcal{N}$ with $[X_i, X_j]=\sum_{k=1}^n c_{ij}^k X_k$, then the structure constants of $\mathcal{N}^{\varphi}$ with respect to $X_1,...,X_n$ are given by $\varphi(c_{ij}^k)$. The $\varphi$-conjugate is well defined up to isomorphism \cite[Corollary 3.6]{Dere}. Let $\overline{\varphi}:\mathcal{N}\to \mathcal{N}^{\varphi}$ be the map $\overline{\varphi}(\sum_{i} a_i X_i)=\sum_i \varphi(a_i)X_i$. Let $N^{\varphi}$ be a simply connected complex Lie group with Lie algebra $\mathcal{N}^{\varphi}$. Let $\hat{\varphi}:N\to N^{\varphi}$ be the map $\hat{\varphi}=\exp_{N^{\varphi}} \circ \overline{\varphi} \circ \exp_N^{-1}$, where $\exp_N$, $\exp_{N^{\varphi}}$ denote the Lie group exponential maps of $N$, $N^\varphi$, respectively. The map $\hat{\varphi}$ is an abstract group isomorphism (see Remark \ref{remark_phi_hat_is_group_isomorphism}), and if $\varphi(c_{ij}^k)=c_{ij}^k$, we have $N=N^{\varphi}$ and $\hat{\varphi}$ is a field automorphism of $N$. For this reason, we call $\hat{\varphi}$ a \textit{field isomorphism} induced by $\varphi\in\textrm{Aut}(\mathbb{C})$.

\begin{remark}
    If $\mathcal{N}$ has a basis whose structure constants $c_{ij}^k$ are in $\mathbb{Z}$, then we have $N=N^\varphi$ for any field automorphism $\varphi$ of $\mathbb{C}$. This is certainly true for semi-simple Lie groups, as well as Iwasawa N-groups. However, this is not true in general (see Example \ref{example_by_Dere}).
\end{remark}

 \begin{definition} \label{definition_of_standard_automorphism}
     Let $N$ be a simply connected complex nilpotent Lie group. An abstract group automorphism $F$ of $N$ is called \textit{standard} if there exists a field automorphism $\varphi$ of $\mathbb{C}$, a central automorphism of $N$, and a complex Lie group isomorphism $\overline{F}:N^\varphi\to N$ such that 
     \begin{align*}
         F=\mu \circ \overline{F} \circ \hat{\varphi},
     \end{align*}
     where $\hat{\varphi}:N\to N^\varphi$ is a field isomorphism induced by $\varphi$.
 \end{definition}

 This definition is related to the base change of affine algebraic groups, and Borel and Tits presented their results that way \cite{Borel_and_Tits}. Here, we do not have to employ the theory of algebraic groups, as discussed above.

 If $N=N^\varphi$ for any $\varphi\in \textrm{Aut}(\mathbb{C})$, then every abstract group automorphism being standard implies partial automatic continuity. However, in general, a standard automorphism may not factor as a product of a complex Lie group automorphism, central automorphism, and field automorphism (see Example \ref{examole_introduction_standard_vs_pac}).

\begin{maintheorem} \label{maintheorem_over_C}
    Let $N$ be a complex simply connected 2-step nilpotent Lie group. Let $\mathcal{N}$ be the Lie algebra of $N$. Let $\mathcal{Z}$ be the center of $\mathcal{N}$. Suppose one of the following conditions holds:
    \begin{enumerate}
        \item $N$ is a generic simply connected complex 2-step nilpotent Lie group with $\dim [N,N] \geq \frac{\dim N-1}{2}$.
        \item There exists a vector space complement $\mathcal{V}$ of $\mathcal{Z}$ in $\mathcal{N}$ and a basis $e_1,...,e_q$ of $\mathcal{V}$ such that $[e_1,. e_j]$, $j\in \{2, ..., q\}$, is linearly independent.
        \item If $\mathcal{N}'$ denotes the Lie algebra of one of the Iwasawa N-groups of a real simple Lie group of rank 1, then $\mathcal{N}=\mathcal{N}'\otimes_\mathbb{R} \mathbb{C}$.
    \end{enumerate}

    Then, every abstract group automorphism of $N$ is standard. In particular, if $\mathcal{N}$ admits a basis whose structure constants lie in $\mathbb{Z}$ (such as the item 3), then $N$ satisfies partial automatic continuity.
\end{maintheorem}

\begin{remark} \label{remark_affine_algebraic_groups}
    In fact, our main technical theorem (Theorem \ref{theorem_concatenation}) is proved over any field $\mathbb{F}$. An interested reader can easily obtain a similar result to Theorem \ref{maintheorem_over_C} for 2-step nilpotent affine algebraic groups over a field of characteristic not 2 using Corollary \ref{corollary_field_of_char_not_2} (cf. Baer correspondence \cite{Baer}), so we only focus on $\mathbb{F}=\mathbb{R},\mathbb{C}$ in the Introduction.
\end{remark}

It is natural to ask:
\begin{question}
    Does there always exist a basis $X_1,...,X_n$ of $\mathcal{N}$ with the structure constants $c_{ij}^k$ such that $\varphi(c_{ij}^k)=c_{ij}^k$, and hence $N=N^\varphi$ and $\hat{\varphi}$ is in fact a field automorphism and $\overline{F}$ is a complex Lie group automorphism?
\end{question}

\begin{example} \label{examole_introduction_standard_vs_pac}
    There exists a standard abstract group automorphism of some simply connected complex 2-step nilpotent Lie group that cannot be decomposed into a central automorphism, a complex Lie group automorphism, and a field automorphism.
\end{example}

Any counterexample to Der\'e's conjecture works (\cite{Demarche}, \cite{Borovoi-deGraaf-Guralnick}); see Propositon \ref{propositon_standard_may_not_be_a_product}. Presumably, generic complex 2-step nilpotent Lie groups may not satisfy partial automatic continuity even if all abstract group automorphisms are standard. The interesting question of when every group automorphism being standard implies partial automatic continuity is not pursued in this paper.

\subsection{Proof Strategy and Some Remarks}
For the proof, we reduce the problem to the Lie algebra level. Let $N$ be a simply connected 2-step nilpotent Lie group with Lie algebra $\mathcal{N}$. Let $\mathcal{Z}$ be the center of $\mathcal{N}$. Take a vector space complement $\mathcal{V}$ so $\mathcal{N}=\mathcal{V}\oplus \mathcal{Z}$. Since $[\mathcal{N}, \mathcal{N}]\subseteq \mathcal{Z}$, there exists a subspace $\mathcal{A}$ such that $\mathcal{Z}=[\mathcal{N}, \mathcal{N}]\oplus \mathcal{A}$. This $\mathcal{A}$ is called an abelian factor of $\mathcal{N}$. 

A bijective additive map that preserves the Lie bracket is called a Lie ring automorphism. Due to Malcev Correspondence (Theorem \ref{correspondence}), an abstract group automorphism $F$ of $N$ corresponds to a Lie ring automorphism $f$ of $\mathcal{N}$. Importantly, $f$ preserves the center and centralizer. 

If $\mathcal{A}=0$, then $\mathcal{N}=\mathcal{V}\oplus \mathcal{Z}$ is a Carnot grading. It follows from Cornulier's much more general result \cite[Corollary 3.6]{Cornulier} that $f$ is a product of a central automorphism $\mu$ and a grading-preserving Lie ring automorphism $\overline{f}$. In our specific setting, the key is that $\overline{f}$ preserves the linear independence (cf. Lemma \ref{lemma_from_arbitrary_to_admissible} (2)). A core result is proved in Section \ref{section_general_results}, and a key lemma is Lemma \ref{lemma_from_arbitrary_to_admissible}. Except for Theorem \ref{isomorphism_theorem_introduction}, all theorems follow straightforwardly from our results in Section \ref{section_general_results}. For Theorem \ref{isomorphism_theorem_introduction}, we use some Lie theory.

In Sections \ref{section_semidirect_product} and \ref{section_abelian_factors}, we reduce the case of $\mathcal{A}\neq 0$ to the case of $\mathcal{A}=0$. 

\vskip\baselineskip
\noindent
\textbf{Remarks about Generalizations.} Section \ref{section_general_results}, \ref{section_semidirect_product}, and \ref{section_abelian_factors} are written for Lie ring automorphisms of finite-dimensional Lie algebras over any field $\mathbb{F}$. See Remark \ref{remark_affine_algebraic_groups}. We emphasize that Theorem \ref{isomorphism_theorem_introduction} does not extend to the affine algebraic group setting; see Example \ref{example_arising_from_field_of_rational_functions}. 

Our proof method does not work for real Lie algebras defined over a $\mathbb{R}$-algebra that is not a field. Given that Question \ref{driving_question} concerns a possible discontinuous automorphism induced by discontinuous ring automorphisms of $\mathbb{R}$-algebras, it is natural to ask if all of our results extend to that setting. This question is not pursued here and remains open.

\subsubsection*{Organization.} In Section \ref{section_automorphisms}, we reduce the problem to the Lie algebra level (Lemma \ref{reduction_to_the_Lie_algebras}), using Malcev correspondence (Theorem \ref{correspondence}). We also describe central, ring, and field automorphisms. Some examples are discussed. Section \ref{section_general_results} establishes Theorem \ref{theorem_concatenation} and Corollary \ref{corollary_concatenation}, which are key tools for our main theorems. Corollary \ref{corollary_concatenation} covers the most general Lie groups treated in this paper.  In Section \ref{section_generic}, Theorem \ref{generic_condition_introduction}, Theorem \ref{maintheorem_criterion}, and the cases (1) and (2) of Theorem \ref{maintheorem_over_C} are proved. In Section \ref{section_counterexamples}, Theorem \ref{pac_up_to_dim_6_introduction} is proved. In Section \ref{section_definitions}, Theorem \ref{temp_theorem1_intro} and the case (3) of Theorem \ref{maintheorem_over_C} are proved. In Section \ref{section_isomorphism_rigidity}, Theorem \ref{isomorphism_theorem_introduction} is proved. In Section \ref{section_semidirect_product}, the semidirect product structure of the group of abstract group automorphisms is explained. In Section \ref{section_abelian_factors}, we reduce the problem to the case without an abelian factor. The results from the last two sections are used in Sections \ref{section_general_results},..., \ref{section_isomorphism_rigidity}, since Sections \ref{section_semidirect_product} and \ref{section_abelian_factors} are independent of other sections.

\subsubsection*{Acknowledments.} I would like to thank my advisor Michael Jablonski for giving helpful suggestions. I thank Jorge Lauret for the questions he asked at the conference ``Symmetry and Geometry in South Florida,'' which led to an improvement of an argument. This work was supported in part by National Science Foundation grant DMS-1906351.

\section{Automorphisms}\label{section_automorphisms}
In this section, we describe a tool we use to study abstract group homomorphisms and automorphisms. Then, central and field automorphisms are explained. We reduce the study to the Lie algebra level.

\subsection{Malcev Correspondence}\label{section_Malcev_correspondence}
A Lie group homomorphism $F$ and its derivative, a Lie algebra homomorphism, $f=dF_e$ relate as $F\circ \exp=\exp \circ f$. For nilpotent Lie groups, Malcev correspondence tells us that a similar relation is still true for a group homomorphism and an additive map preserving the Lie bracket. In other words, forgetting a smooth structure at the Lie group level corresponds to forgetting $\mathbb{R}$-linearity at the Lie algebra level. Before we state the result, we make the following definition.

\begin{definition}
Let $\mathfrak{g}, \mathfrak{h}$ be finite-dimensional nilpotent Lie algebras over $\mathbb{R}$.

A map $f:\mathfrak{g}\to \mathfrak{h}$ is called a \textit{Lie ring homomorphism} if for any $x,y \in \mathfrak{g}$, $f(x+y)=f(x)+f(y)$ and $f([x,y])=[f(x),f(y)]$, that is, $f$ is additive and $f$ preserves the Lie bracket. A bijective Lie ring homomorphism is called a \textit{Lie ring isomorphism}. When $\mathfrak{g}=\mathfrak{h}$, a Lie ring isomorphism is called a \textit{Lie ring automorphism}.
\end{definition}

A Lie ring is defined to be a Lie algebra over $\mathbb{Z}$, and any Lie algebra can be thought of as a Lie ring by forgetting the $\mathbb{R}$-linear structure. The name Lie ring homomorphisms comes from this fact. Note, a Lie ring homomorphism is in fact $\mathbb{Q}$-linear.

Our main tool is a generalized version of Malcev correspondence (\cite{Malcev_Nilpotent_Torsion_Free_Groups}):

\begin{theorem}\label{correspondence}
Let $\mathfrak{g}, \mathfrak{h}$ be finite-dimensional nilpotent Lie algebras over $\mathbb{R}$. Let $G, H$ be simply connected nilpotent Lie groups corresponding to $\mathfrak{g}, \mathfrak{h}$, respectively.

Then, for any abstract group homomorphism $F:G\to H$, the map $\mathfrak{g}\to \mathfrak{h}$ defined by 
\begin{equation*}
    f=\exp_H^{-1}\circ F\circ \exp_G
\end{equation*}
is a Lie ring homomorphism. For any Lie ring homomorphism $f:\mathfrak{g}\to \mathfrak{h}$, the map $G\to H$ defined by
\begin{align*}
   F= \exp_H\circ f\circ \exp_G^{-1}
\end{align*}
is a group homomorphism. In particular, if $G=H$, then the group of group automorphisms and the group of Lie ring automorphisms are isomorphic as a group via $f\to F=\exp_G\circ f\circ \exp_G^{-1}$.
\end{theorem}

\begin{proof}
    Since the exponential map is a diffeomorphism, simply connected nilpotent Lie groups $G, H$ are Malcev groups, as defined in \cite[p.98]{Amayo_Stewart}. The associated Lie algebras $\mathscr{L}(G), \mathscr{L}(H)$ over $\mathbb{Q}$ defined in \cite[p.103]{Amayo_Stewart} are isomorphic to the Lie algebras $\mathfrak{g}, \mathfrak{h}$ regarded as the Lie algebras over $\mathbb{Q}$, respectively. By \cite[Theorem 5.2(d), p.108]{Amayo_Stewart}, the result follows.  
\end{proof}
A more detailed proof is also found at \cite[Chapter 6, Proposition 4]{Segal}. In \cite{Segal}, the setting is different, but their proof works for our setting. 

Since the exponential map is a diffeomorphism, we have the following.

\begin{lemma}\label{continuity_of_F_and_f}
    Let $N$ be a simply connected nilpotent Lie group. Let $\mathcal{N}$ be its Lie algebra. Let $F:N\to N$ be a group homomorphism. Let $f=\exp^{-1}\circ F\circ \exp:\mathcal{N}\to \mathcal{N}$ be the corresponding Lie ring homomorphism from Theorem \ref{correspondence}.

    Then, $F$ is continuous if and only if $f$ is continuous. If $F$ is continuous, then $F$ is a Lie group homomorphism. If $f$ is continuous, then $f$ is a Lie algebra homomorphism. Hence, $F$ is a Lie group homomorphism if and only if $f$ is a Lie algebra homomorphism.
\end{lemma}

\begin{proof}
    The first statement holds because the exponential map is a diffeomorphism. A continuous group homomorphism is smooth. Suppose $f$ is continuous. Then, since $f$ is $\mathbb{Q}$-linear, $f$ must be $\mathbb{R}$-linear.
\end{proof}

\subsection{Central Automorphisms}\label{section_central_automorphisms}

Given a group $N$, an automorphism $F$ is called a \textit{central automorphism} if for any $x\in N$, $x^{-1}F(x)\in Z$, the center of $N$. Central automorphisms are precisely the centralizers of the group of inner automorphisms. The group of central automorphisms forms a normal subgroup of the group of abstract group automorphisms. This definition coincides with definitions in \cite{Gibbs} and \cite{Khor}. Let $N$ be a simply connected nilpotent Lie group with center $Z$. Let $\mathcal{N}$ be its Lie algebra with center $\mathfrak{z}$. We have the following equivalence.

\begin{proposition}\label{characterization_of_central_automorphism}
     Let $F$ be a group automorphism. Let $f=\exp^{-1}\circ F\circ \exp$ be the corresponding Lie ring automorphism from Theorem \ref{correspondence}.

    Then, $F$ is a central automorphism if and only if $f(x)-x\in \mathfrak{z}$ for any $x\in \mathcal{N}$.

\end{proposition}

\begin{proof}
    Recall that $\exp$ is a diffeomorphism and $\exp(\mathfrak{z})=Z$. Suppose that $F$ is central. Let $\lambda(x)=x^{-1}F(x)$. Then, $F(x)=x\lambda(x)$. Since $\lambda(x) \in Z$, $\exp^{-1}\circ F(x)=\exp^{-1}(x)+\exp^{-1}\circ \lambda(x)$. Thus, $f(y)=y+\exp^{-1}\circ\lambda \circ \exp(y)$ for $y=\exp^{-1}(x)$. The converse is similar.
\end{proof}

\begin{definition}
    Let us call a Lie ring automorphism $f$ a \textit{central automorphism} if $f(x)-x\in \mathfrak{z}$ for any $x\in \mathcal{N}$.
\end{definition}

Take a subspace $\mathcal{V}$ of $\mathcal{N}$ such that $\mathcal{N}=\mathcal{V}\oplus [\mathcal{N}, \mathcal{N}]$. A central automorphism $f$ corresponds to an additive map from $\mathcal{V}$ to $\mathfrak{z}$, as explained below.
 
Let $f$ be a central Lie ring automorphism. Let $\mu(x)=f(x)-x$ so that $f=\textrm{id}_{\mathcal{N}}+\mu$. Then, $\mu(\mathcal{V})\subseteq \mathfrak{z}$. The central part $\mu(x)=f(x)-x$ is additive and 0 on the commutator ideal. Indeed, $[x,y]+\mu([x,y])=f([x,y])=[x+\mu(x), y+\mu(y)]=[x,y]$. Thus, $\mu$ can be thought of as an additive map from $\mathcal{V}$ to the center $\mathfrak{z}$. Conversely,

\begin{lemma}\label{characterization_central_automorphism}
    Let $\mathcal{N}$ be a nilpotent Lie algebra with center $\mathfrak{z}$.
    Let $\mu:\mathcal{N}\to \mathfrak{z}$ be an additive map such that $\mu|_{[\mathcal{N},\mathcal{N}]}=0$. If $\mu^2=0$, then $f=\textrm{id}_{\mathcal{N}}+\mu$ is a central automorphism with the inverse $\textrm{id}_{\mathcal{N}}-\mu$.
\end{lemma}

\begin{proof}
     Since $f$ acts on the commutator ideal trivially and the image of $\mu$ is in the center, we have $f([x,y])=[x,y]=[f(x),f(y)]$. By design, $f(x)-x=\mu(x)\in \mathfrak{z}$ for any $x\in \mathcal{N}$. Since $\mu^2=0$, $f\circ(\textrm{id}_{\mathcal{N}}-\mu)=\textrm{id}_{\mathcal{N}}$ and $(\textrm{id}_{\mathcal{N}}-\mu)\circ f=\textrm{id}_{\mathcal{N}}$.  
\end{proof} 

The map $\mu$ is a $\mathbb{Q}$-linear map. Under the axiom of choice, there exists a $\mathbb{Q}$-basis of $\mathcal{V}$ and $\mu$ rearranges such a basis. Thus, they are typically discontinuous. See \cite[Corollary 5.2.2]{Kuczma} for more details.

\begin{proposition} \label{discontinuous_automorphisms_for_nilpotent_Lie_groups}
    Let $N$ be a simply connected non-abelian nilpotent Lie group, and let $\mathcal{N}$ be its Lie algebra. Let $\mathcal{V}$ be a subspace such that $\mathcal{N}=\mathcal{V}\oplus [\mathcal{N}, \mathcal{N}]$. Since $\mathcal{N}$ is nilpotent, $\mathcal{V}$ is not zero. Let $\mathfrak{z}$ be the center of $\mathcal{N}$. Since $\mathcal{N}$ is nilpotent, $\mathfrak{z}\cap [\mathcal{N}, \mathcal{N}]$ is not zero. 
    Let $\mu:\mathcal{V}\to \mathfrak{z} \cap [\mathcal{N}, \mathcal{N}]$ be a discontinuous additive map. Define a map $f:\mathcal{N}\to \mathcal{N}$ by $f(X+Y)=X+Y+\mu(X)$ for $X\in \mathcal{V}$ and $Y\in [\mathcal{N}, \mathcal{N}]$. Then, $f$ is a discontinuous Lie ring automorphism, and $F=\exp \circ f \circ \exp^{-1}:N\to N$ is a discontinuous central automorphism.
\end{proposition}

\begin{proof}
    Extend $\mu$ to a $\mathbb{Q}$-linear map from $\mathcal{N}$ to $\mathfrak{z}\cap [\mathcal{N}, \mathcal{N}]$ by setting $\mu|_{[\mathcal{N},\mathcal{N}]}=0$. Since $\mu^2=0$, it follows from Lemma \ref{characterization_central_automorphism} that $f$ is a central automorphism. By Proposition \ref{characterization_of_central_automorphism}, $F$ is a central automorphism. Since $\mu$ is discontinuous, so is $f$. It follows from Theorem \ref{correspondence} and Lemma \ref{continuity_of_F_and_f} that $F=\exp \circ f \circ \exp^{-1}$ is a discontinuous group automorphism.
\end{proof}

\subsection{Ring and Field Automorphisms}\label{section_field_automorphisms}
In Question \ref{driving_question}, we consider ring automorphisms of simply connected nilpotent Lie groups. We formally define them as follows.

Let $\mathcal{N}$ be a real nilpotent Lie algebra. Suppose there exists a $\mathbb{R}$-algebra $A$ such that $\mathcal{N}$ is defined over $A$, i.e., there exists an $A$-basis $\{e_i\}_{i=1}^n$ of $\mathcal{N}$. Let $\varphi$ be a ring automorphism of $A$. Define $\overline{\varphi}:\mathcal{N}\to \mathcal{N}$ by 
\begin{align*}
\overline{\varphi}\left( \sum_{i=1}^n x_i e_i \right) = \sum_{i=1}^n \varphi(x_i) e_i.
\end{align*}
Suppose that $\varphi$ fixes the structure constants. Then, $\overline{\varphi}$ is a Lie ring automorphism of the real Lie algebra $\mathcal{N}$. We call $\overline{\varphi}$ a \textit{ring automorphism}. When $A$ is the field $\mathbb{C}$, we call $\overline{\varphi}$ a \textit{field automorphism}

Let $N$ be a simply connected nilpotent Lie group with its Lie algebra $\mathcal{N}$. By Theorem \ref{correspondence}, $\exp\circ \overline{\varphi}\circ \exp^{-1}$ is a group automorphism if $\varphi$ is a ring automorphism. We call it a \textit{ring automorphism} as well. When $\mathbb{F}$ is a field $\mathbb{C}$, it is called a \textit{field automorphism}.

\begin{remark}
    The construction depends on the basis. If $\mathcal{N}$ admits a basis with all the structure constants in $\mathbb{Z}$, we can use such a basis.
\end{remark}

\begin{remark}
    This definition generalizes the field automorphisms defined in \cite{Gibbs} and \cite{Khor}. 
\end{remark}

\subsection{Lie Algebras over a Field}
The main technical result (Theorem \ref{theorem_concatenation}) holds over any field $\mathbb{F}$. For this reason, let $\mathbb{F}$ be any field, and we will discuss some properties and examples of Lie algebras over $\mathbb{F}$. Any Lie algebra $\mathcal{N}$ is naturally a Lie ring. 

\begin{remark}
    Ring automorphisms and central automorphisms are defined in the same way as in the case of $\mathbb{F}=\mathbb{R}$.
\end{remark}

First, we record the following useful lemma.
\begin{lemma} \label{center_goes_to_center}
    Let $\mathfrak{g}, \mathfrak{h}$ be Lie rings, and let $\mathfrak{z}(\mathfrak{g}), \mathfrak{z}(\mathfrak{h})$ be the center of $\mathfrak{g}, \mathfrak{h}$, respectively.

    If $f:\mathfrak{g}\to \mathfrak{h}$ is a Lie ring isomorphism, that is, a bijective additive map preserving the Lie bracket, then $f(\mathfrak{z}(\mathfrak{g}))=\mathfrak{z}(\mathfrak{h})$.
\end{lemma}
\begin{proof}
    Let $x\in \mathfrak{z}(\mathfrak{g})$. Then, for any $y\in \mathfrak{h}$, $[f(x),y]=[f(x),f(f^{-1}(y))]=f([x,f^{-1}(y)])=f(0)=0$. Thus, $f(x)\in \mathfrak{z}(\mathfrak{h})$. Thus, $f(\mathfrak{z}(\mathfrak{g}))\subseteq \mathfrak{z}(\mathfrak{h})$. By reversing the role, one also has $f^{-1}(\mathfrak{z}(\mathfrak{h}))\subseteq \mathfrak{z}(\mathfrak{g})$.
\end{proof}

Our framework is over a field, so we focus on field automorphisms. A related notion to field automorphisms is semilinearity. Recall the following definition.

\begin{definition}
    Let $\mathbb{F}$ be any field, and $V$ be a vector space over $\mathbb{F}$. Let $f:V\to V$ be an additive map. Suppose that there exists a field homomorphism $\varphi$ of $\mathbb{F}$ such that for any $t\in \mathbb{F}$ and $x\in V$, $f(tx)=\varphi(t)f(x)$.

    Then, the map $f$ is called $\mathbb{F}$\textit{-semilinear} or $\varphi$\textit{-linear}.
\end{definition}

This is a classical notion commonly used in affine and projective geometry. Now suppose $f$ and $\varphi$ are bijective. It is clear that for a fixed basis $\{e_i\}_{i=1}^n$ of $V$, $\varphi$ induces a $\varphi$-linear map \begin{align*}
    \overline{\varphi}\left( \sum_{i=1}^n x_i e_i \right) = \sum_{i=1}^n \varphi(x_i) e_i.
\end{align*}
Set $\overline{f}:=f\circ \overline{\psi}$, where $\psi=\varphi^{-1}$. Then $\overline{f}$ is $\mathbb{F}$-linear. Indeed, for any $t\in \mathbb{F}$ and $j=1,...,n$, we have $\overline{f}(te_j)=f(\psi(t)e_j)=\varphi(\psi(t))f(e_j)=tf(e_j)=t\overline{f}(e_j)$. By using the additivity, $\overline{f}$ is $\mathbb{F}$-linear, and $f=\overline{f}\circ \overline{\varphi}$. 

Recall that we defined the $\varphi$-conjugate $\mathcal{N}^{\varphi}$ of a complex Lie algebra $\mathcal{N}$ for a field automorphism $\varphi\in \textrm{Aut}(\mathbb{C})$ in Section \ref{section_introduction_complex_case}. For a Lie algebra over $\mathbb{F}$ and $\varphi\in \textrm{Aut}(\mathbb{F})$, $\mathcal{N}^\varphi$ is defined in the same way (see \cite[Section 3]{Dere}). For a fixed basis $\{e_i\}_{i=1}^n$ of $\mathcal{N}$, the map $\overline{\varphi}:\mathcal{N}\to \mathcal{N}^{\varphi}$ defined by $\overline{\varphi}(\sum_{i=1}^n x_ie_i)=\sum_{i=1}^n \varphi(x_i)e_i$ is called a \textit{field isomorphism} induced by $\varphi$. (Recall $\mathcal{N}=\mathcal{N}^{\varphi}$ as vector spaces.)

\begin{remark} \label{remark_phi_hat_is_group_isomorphism}
    The map $\overline{\varphi}$ is a Lie ring isomorphism by design. When $\mathbb{F}=\mathbb{C}$, $\hat{\varphi}=\exp_{N^\varphi}\circ \overline{\varphi}\circ \exp_N^{-1}$, as defined in Section \ref{section_introduction_complex_case}, is an abstract group isomorphism due to Theorem \ref{correspondence}.
\end{remark}

\begin{proposition} \label{existence_of_semilinear_Lie_ring_automorphism}
    Let $\mathcal{N}$ be a Lie algebra over $\mathbb{F}$. Let $\varphi$ be a field automorphism of $\mathbb{F}$. Then,
    \begin{enumerate}
        \item $\mathcal{N}$ admits a $\varphi$-linear Lie ring automorphism if and only if $\mathcal{N}$ is isomorphic to $\mathcal{N}^{\varphi}$ as Lie algebras over $\mathbb{F}$.
        \item Any $\varphi$-linear Lie ring automorphism $f$ is of the form $f=\overline{f}\circ \overline{\varphi}$, where $\overline{f}$ is a Lie algebra isomorphism $\overline{f}:\mathcal{N}^{\varphi}\to \mathcal{N}$ and $\overline{\varphi}:\mathcal{N}\to \mathcal{N}^\varphi$ is a field isomorphism.
        \item If there exists a basis whose structure constants are fixed by $\varphi$, then any $\varphi$-linear Lie ring automorphism $f$ is of the form $f=\overline{f}\circ \overline{\varphi}$, where $\overline{f}$ is a Lie algebra automorphism of $\mathcal{N}$ and $\overline{\varphi}$ is a field automorphism of $\mathcal{N}$.
    \end{enumerate}

\end{proposition}
\begin{proof}
     Let $e_1,..., e_n$ be a fixed basis of $\mathcal{N}$. Define $\overline{\varphi}$ by $\overline{\varphi}\left( \sum_{i=1}^n x_i e_i \right) = \sum_{i=1}^n \varphi(x_i) e_i$. Then, $\overline{\varphi}:\mathcal{N}\to \mathcal{N}^{\varphi}$ is a Lie ring isomorphism by design. 
     
     First assume $f$ is a $\varphi$-linear Lie ring automorphism and show the item 2. Let $\overline{f}=f\circ \overline{\varphi}^{-1}:\mathcal{N}^{\varphi}\to \mathcal{N}$. Then, $\overline{f}$ is a Lie ring isomorphism as well. Moreover, $\overline{f}$ is $\mathbb{F}$-linear by the same discussion as above. Thus, $\overline{f}$ is a Lie algebra isomorphism. This shows the item 2 and the forward direction of the item 1.

    To show the other direction of the item 1, assume there is a Lie algebra isomorphism $\overline{f}:\mathcal{N}^\varphi\to \mathcal{N}$. Then set $f=\overline{f} \circ \overline{\varphi}:\mathcal{N}\to \mathcal{N}$. The map $f$ is the desired $\varphi$-linear Lie ring automorphism.

    The item 3 immediately follows from the item 2 because for such a basis, $\mathcal{N}=\mathcal{N}^\varphi$.
\end{proof}

Let $\mathcal{N}$ be a complex Lie algebra and $\sigma\in \textrm{Aut}(\mathbb{C})$ be the complex conjugation. Der\'e (\cite[Conjecture 1]{Dere}) conjectured that if $\mathcal{N}$ is isomorphic to $\mathcal{N}^{\sigma}$, then $\mathcal{N}$ is defined over $\mathbb{R}$. Recently, Demarche \cite{Demarche} proved the existence of a 10-dimensional complex 2-step nilpotent counterexample, and an explicit counterexample was found by Borovoi--de Graaf--Guralnick \cite{Borovoi-deGraaf-Guralnick}. 

\begin{proposition} \label{propositon_standard_may_not_be_a_product}
    Let $\mathcal{N}$ be a complex 2-step nilpotent Lie algebra that is a counterexample to \cite[Conjecture 1]{Dere} (\cite{Demarche}, \cite{Borovoi-deGraaf-Guralnick}). Let $N$ be a simply connected complex 2-step nilpotent Lie group with Lie algebra $\mathcal{N}$.
    
    Then, $N$ admits a standard group automorphism that is not a product of a central automorphism, complex Lie group automorphism, and field automorphism.
\end{proposition}

\begin{proof}
    Since $\mathcal{N}$ is isomorphic to $\mathcal{N}^{\sigma}$, $\mathcal{N}$ admits a $\sigma$-linear Lie ring automorphism $f$ by Proposition \ref{existence_of_semilinear_Lie_ring_automorphism}. On the other hand, if $\sigma$ induces a field automorphism $\overline{\sigma}$ of $\mathcal{N}$, then by definition, there exists a basis of $\mathcal{N}$ whose structure constants $c_{ij}^k$ are fixed by $\sigma$. Since $\sigma(c_{ij}^k)=c_{ij}^k$, we have $c_{ij}^k\in \mathbb{R}$, i.e., $\mathcal{N}$ is defined over $\mathbb{R}$. This is a contradiction. Thus, $\sigma$ cannot induce a field automorphism of $\mathcal{N}$.

    Suppose $f=\mu \circ \overline{f}\circ \overline{\varphi}$, where $\mu$ is a central automorphism, $\overline{f}$ is a complex Lie algebra automorphism, and $\overline{\varphi}$ is a field automorphism of $\mathcal{N}$ induced by some $\varphi\in \textrm{Aut}(\mathbb{C})$. As discussed, $\varphi\neq \sigma$. Observe, all three maps descend to $\mathcal{N}/\mathcal{Z} \to \mathcal{N}/\mathcal{Z}$, where $\mathcal{Z}$ is the center (Lemma \ref{center_goes_to_center}). Since $\mu(x)-x\in \mathcal{Z}$, $\mu$ descends to the identity map, so $f=\overline{f}\circ \overline{\varphi}:\mathcal{N}/\mathcal{Z}\to \mathcal{N}/\mathcal{Z}$. Here, we abuse the notation and use $f, \overline{f}, \overline{\varphi}$ for the descended maps.
    
    For any $t\in \mathbb{R}$ and $x\in \mathcal{N}$, $f(tx)=\sigma(t)f(x)=tf(x)$, so $f$ is $\mathbb{R}$-linear, and hence continuous. Thus, $f$ and $\overline{f}$ are both smooth, and $\overline{\varphi}:\mathcal{N}/\mathcal{Z}\to \mathcal{N}/\mathcal{Z}$ has to be smooth. This means $\varphi$ is continuous, so $\varphi=\textrm{Id}$ or $\sigma$, but $\varphi\neq \sigma$. Thus, $\varphi=\textrm{Id}$, and $f=\overline{f}$ on $\mathcal{N}/\mathcal{Z}$. However, for any $x\in \mathcal{N}\setminus\mathcal{Z}$, $f(\sqrt{-1}(x+\mathcal{Z}))=-\sqrt{-1}f(x+\mathcal{Z})$, while $\overline{f}$ is $\mathbb{C}$-linear. This is a contradiction. Thus, $F=\exp \circ f\circ \exp^{-1}$ is the desired standard automorphism (see Lemma \ref{reduction_to_the_Lie_algebras}).
\end{proof}

Note that $\mathcal{N}$ is not necessarily isomorphic to $\mathcal{N}^{\varphi}$ for all field automorphisms $\varphi$ of $\mathbb{F}$. The following example by J. Der\'e is interesting.
\begin{example}(\cite[Example 5.1, Proposition 5.2]{Dere}) \label{example_by_Dere}
Let $\sigma(z)=\overline{z}$ be the complex conjugation.
    There is a 10-dimensional 2-step nilpotent Lie algebra $\mathcal{N}$ such that $\mathcal{N}$ is not isomorphic to $\mathcal{N}^{\sigma}$. 
\end{example}

Thus, even in the 2-step nilpotent Lie algebras over $\mathbb{C}$, the complex conjugation $\sigma$ does not always give rise to a $\sigma$-linear Lie ring automorphism. This issue arises as nilpotent Lie algebras over $\mathbb{C}$ do not necessarily admit nice bases ($\mathbb{Z}$-basis, $\mathbb{R}$-basis, etc.), unlike semisimple Lie algebras.

Weirdly, for some field $\mathbb{F}$, $\textrm{heis}_3(\mathbb{F})$ can be isomorphic to a Lie algebra $\tilde{\mathcal{N}}$ over $\mathbb{F}$ with $\dim_\mathbb{F} \tilde{\mathcal{N}}=6$ as Lie rings, showing that Theorem \ref{isomorphism_theorem_introduction} does not extend to affine algebraic groups over such a field.
Let $\mathbb{K}$ be a field. Let $\mathbb{F}=\mathbb{K}(x)$ be the field of rational functions over $\mathbb{K}$. Let $\mathbb{E}=\mathbb{K}(x^2)$. The key is that there exists a field homomorphism $\varphi:\mathbb{F}\to \mathbb{F}$ such that $\varphi(1)=1$, $\varphi(x)=x^2$. The image of $\varphi$ is $\mathbb{E}$. Thus $\mathbb{F}=\mathbb{E}\oplus x\mathbb{E}$, where $\mathbb{E}\cong \mathbb{F}$ as fields (and hence as $\mathbb{Z}$-modules). Hence, $\mathbb{F}$ is somewhat naturally isomorphic to $\mathbb{F}^2$ as $\mathbb{Z}$-modules. This carries over to $\textrm{heis}_3(\mathbb{F})$.

\begin{example} \label{example_arising_from_field_of_rational_functions}
    Let $\mathcal{N}=\mathbb{F}\textrm{-span}\{e_1, e_2, Z_1\}$ be the 3-dimensional Heisenberg algebra over $\mathbb{F}$ with $[e_1,e_2]=Z_1$. Let $\tilde{\mathcal{N}}=\mathbb{F}\textrm{-span}\{\tilde{e}_1, \tilde{e}_2, \tilde{e}_3, \tilde{e}_4, \tilde{Z}_1, \tilde{Z}_2\}$ be a 6-dimensional 2-step nilpotent Lie algebra over $\mathbb{F}$ with nontrivial Lie brackets $[\tilde{e}_1, \tilde{e}_2]=\tilde{Z}_1, [\tilde{e}_1, \tilde{e}_3]=0, [\tilde{e}_1, \tilde{e}_4]=\tilde{Z}_2, [\tilde{e}_2, \tilde{e}_3]=-\tilde{Z}_2, [\tilde{e}_2, \tilde{e}_4]=0, [\tilde{e}_3, \tilde{e}_4]=x\tilde{Z}_1$. Then, $\tilde{\mathcal{N}}$ and $\mathcal{N}$ are isomorphic as Lie rings, although they are not isomorphic as Lie algebras over $\mathbb{F}$.

    First, we construct a bijective map $f:\tilde{\mathcal{N}}\to \mathcal{N}$. Let $\varphi:\mathbb{F}\to \mathbb{F}$ be a field homomorphism such that $\varphi(1)=1$ and $\varphi(x)=x^2$. Then, define an additive map $f:\tilde{\mathcal{N}}\to \mathcal{N}$ by $f(t\tilde{e}_1)=\varphi(t)e_1, f(t\tilde{e}_2)=\varphi(t)e_2, f(t\tilde{e}_3)=x\varphi(t)e_1, f(t\tilde{e}_4)=x\varphi(t)e_2, f(t\tilde{Z}_1)=\varphi(t)Z_1$, and $f(t\tilde{Z}_2)=x\varphi(t)Z_1$. The key is that $\textrm{Im}\varphi = \mathbb{K}(x^2)$, and $[\mathbb{F}:\mathbb{K}(x^2)]=2$. In particular, $1, x$ are $\mathbb{K}(x^2)$-basis of $\mathbb{F}$. Thus, $f$ is bijective. 

    Next, we show $f$ is a Lie ring homomorphism. Since $f$ is additive, it suffices to show that $f([t\tilde{e}_i, s\tilde{e}_j])=[f(t\tilde{e}_i), f(s\tilde{e}_j)]$ for any $i\leq j$. This is a straightforward computation. The key is $f([t\tilde{e}_3, s\tilde{e}_4])=f(tsx\tilde{Z}_1)=\varphi(tsx)Z_1=\varphi(t)\varphi(s)x^2 Z_1$, and $[f(t\tilde{e}_3), f(s\tilde{e}_4)]=[x\varphi(t)e_1, x\varphi(s)e_2]=x^2\varphi(t)\varphi(s)Z_1$. This shows that $f$ is a Lie ring isomorphism. On the other hand, $\dim_{\mathbb{F}}\tilde{\mathcal{N}}=6$, while $\dim_{\mathbb{F}}\mathcal{N}=3$. Thus, they are not isomorphic as Lie algebras over $\mathbb{F}$. If $\textrm{char}(\mathbb{F})=\textrm{char}(\mathbb{K})\neq 2$, there exist affine algebraic groups $\tilde{N}, N$ associated to $\tilde{\mathcal{N}}, \mathcal{N}$; see the description below Corollary \ref{corollary_concatenation}. The groups $\tilde{N}, N$ are abstractly isomorphic, but not isomorphic as algebraic groups (see Remark \ref{remark_Baer_correspondence}).
\end{example}

\begin{remark}
    The Lie algebra $\tilde{\mathcal{N}}$ is $L_{6,22}(\epsilon)$ with $\epsilon=x$ in the classification (\cite{Graaf_char_not_2}, \cite{Cicalo_Serena_Graaf}). Note, $L_{6,22}(\epsilon)$ arises as a family. 
\end{remark}

Combining the results of this section, we reduce our study to the Lie algebra level.

\begin{definition}\label{def_of_partial_automatic_continuity_at_the_Lie_algebra_level}
We say that \textit{a real nilpotent Lie algebra} $\mathcal{N}$ \textit{satisfies partial automatic continuity} if any Lie ring automorphism is a product of a central automorphism, Lie algebra automorphism, and ring automorphism.
\end{definition}

\begin{lemma}\label{reduction_to_the_Lie_algebras}
    Let $\mathbb{F}=\mathbb{R}$ or $\mathbb{C}$. Let $\mathcal{N}$ be a nilpotent Lie algebra over $\mathbb{F}$, and let $N$ be the simply connected nilpotent Lie group with Lie algebra $\mathcal{N}$. Let $F:N\to N$ be a group automorphism and let $f=\exp^{-1}\circ F\circ \exp:\mathcal{N}\to \mathcal{N}$ be the corresponding Lie ring automorphism (Theorem \ref{correspondence}). 
    
    Then, 
    \begin{itemize}
        \item The group $N$ satisfies partial automatic continuity if and only if $\mathcal{N}$ satisfies partial automatic continuity.
        \item The map $F$ is a product of a central automorphism and a Lie group automorphism if and only if $f$ is a product of a central automorphism and a Lie algebra automorphism.
        \item Let $\mathbb{F}=\mathbb{C}$. The map $F$ is standard if and only if $f$ is a product of a central automorphism and a $\varphi$-linear Lie ring automorphism for some $\varphi\in \textrm{Aut}(\mathbb{F})$.
    \end{itemize}

\end{lemma}

\begin{proof}[Proof of Lemma \ref{reduction_to_the_Lie_algebras}]
    Let $\mu:N\to N$ be a central automorphism. Let $\overline{F}:N\to N$ be a Lie group automorphism. By definition, a ring automorphism $\hat{\varphi}$ of $N$ has the form $\hat{\varphi}=\exp\circ \overline{\varphi}\circ \exp^{-1}$ as defined above. Then, if $F= \mu \circ \overline{F}\circ \hat{\varphi}$, we have $f=\exp^{-1}\circ \mu \circ \overline{F}\circ \hat{\varphi} \circ \exp=(\exp^{-1}\circ \mu \circ \exp) \circ (\exp^{-1}\circ \overline{F}\circ \exp)\circ \overline{\varphi}$. By Lemma \ref{continuity_of_F_and_f} and Proposition \ref{characterization_of_central_automorphism}, $f$ is a product of a central automorphism, Lie algebra automorphism, and ring automorphism. The converse is similar. This proves the items 1 and 2. Similarly, the item 3 follows from Proposition \ref{existence_of_semilinear_Lie_ring_automorphism} (2) and Definition \ref{definition_of_standard_automorphism}. 
\end{proof}

\section{General Results}\label{section_general_results}

\subsection{Technical Lemmas}
\begin{definition} \label{definition_admissible_Lie_algebras}
    Let $\mathbb{F}$ be a field. Let $\mathcal{N}$ be a 2-step nilpotent Lie algebra over $\mathbb{F}$ with center $\mathcal{Z}$. Then, we say $\mathcal{N}$ is \textit{admissible relative to $\{e_1,...,e_q\}$} if there is an $\mathbb{F}$-subspace $\mathcal{V}$ of $\mathcal{N}$ and an $\mathbb{F}$-basis $e_1,...,e_q$ of $\mathcal{V}$ such that $\mathcal{N}=\mathcal{V}\oplus \mathcal{Z}$ as an $\mathbb{F}$-vector space, and $[e_1,e_j], j\in \{2,...,q\}$ is linearly independent. 

    We say $\mathcal{N}$ is \textit{admissible} if it is admissible relative to some $\{e_1, ..., e_q\}$ as above. 
\end{definition}

We first state a technical lemma about a Lie ring homomorphism from an arbitrary 2-step nilpotent Lie algebra $\tilde{\mathcal{N}}$ into an admissible 2-step nilpotent Lie algebra $\mathcal{N}$. 

\begin{setting}\label{setting_technical_lemma}
    Let $\mathbb{F}$ be a field. Let $\mathcal{N}$ be a 2-step nilpotent Lie algebra over $\mathbb{F}$. Let $\mathcal{Z}$ be the center of $\mathcal{N}$. Suppose that $\mathcal{N}$ is admissible relative to $\{e_1,...,e_q\}$. Let $\tilde{\mathcal{N}}$ be an arbitrary 2-step nilpotent Lie algebra over $\mathbb{F}$. Let $f:\tilde{\mathcal{N}}\to \mathcal{N}$ be a Lie ring homomorphism. Suppose that for each $j\in \{1,...,q\}$, there exists some $\tilde{e}_j\in \tilde{\mathcal{N}}$ such that $f(\tilde{e}_j)=e_j$. 
\end{setting} 
\begin{remark}
    A surjective Lie ring homomorphism satisfies the condition above. For a technical reason, we do not want to assume that $f$ is surjective.
\end{remark}

Note, $\tilde{e}_1, ..., \tilde{e}_q$ are just non-zero elements of $\tilde{\mathcal{N}}$ whose image under $f$ are $e_1,...,e_q$. For $X\in \tilde{\mathcal{N}}$, write $f(X)=\sum_{j=1}^q f_{e_j}(X)e_j+f_{\mathcal{Z}}(X)$, where $f_{e_j}(X)\in \mathbb{F}$ and $f_{\mathcal{Z}}(X)\in \mathcal{Z}$. 

\begin{lemma} \label{lemma_from_arbitrary_to_admissible}
    Under Setting \ref{setting_technical_lemma}, let $\tilde{\mathcal{V}}$ be the subspace $\mathbb{F}$-span$\{\tilde{e_1},...,\tilde{e_q}\}$. Let $\varphi(t)=f_{e_1}(t\tilde{e}_1)$. Then,
    \begin{enumerate}
        \item $\varphi$ is an injective field homomorphism of $\mathbb{F}$, and  $f(t\tilde{e}_j)\in \varphi(t)e_j + \mathcal{Z}$ for all $j=1,...,q$.
        \item $[\tilde{e}_1, \tilde{e}_j], j\in \{2,...,q\}$ is linearly independent, and $\tilde{e}_1,...,\tilde{e}_q$ is linearly independent.
        \item For any $x\in \tilde{\mathcal{V}}$ and $t\in \mathbb{F}$, $f(tx)\in \varphi(t)f(x)+\mathcal{Z}$.
    \end{enumerate}
\end{lemma} 
\begin{remark}
    This lemma implies $f$ is $\varphi$-linear up to center when it is restricted to a possibly smaller Lie subalgebra generated by $\tilde{\mathcal{V}}$, i.e., $\tilde{\mathcal{V}}+ [\tilde{\mathcal{V}},\tilde{\mathcal{V}}]$. 
\end{remark}
We give a proof of Lemma \ref{lemma_from_arbitrary_to_admissible} here. The main ingredients
consist of four sublemmas. In all of these sublemmas, we let $\mathcal{N}, \tilde{\mathcal{N}}$, and $f:\tilde{\mathcal{N}}\to \mathcal{N}$ satisfy the hypothesis in Setting \ref{setting_technical_lemma} and we use the notations above.

\begin{sublemma}\label{sublemma1_for_lemma_from_arbitrary_to_admissible}
    For any $t\in \mathbb{F}$, $f(t\tilde{e}_1)\in \varphi(t)e_1 + \mathcal{Z}$.
\end{sublemma}
\begin{proof}
    We have $0=f([\tilde{e}_1, t\tilde{e}_1])=[e_1, f(t\tilde{e}_1)]=\sum_{j=2}^q f_{e_j}(t\tilde{e}_1)[e_1,e_j]$ as $f(t\tilde{e}_1)=\sum_{j=1}^q f_{e_j}(t\tilde{e}_1)e_j+f_{\mathcal{Z}}(t\tilde{e}_1)$ and $f_{\mathcal{Z}}(t\tilde{e}_1)$ is central. By assumption, $[e_1, e_2],...,[e_1,e_q]$ is linearly independent. Thus, we have $f_{e_j}(t\tilde{e}_1)=0$ for all $t\in \mathbb{F}$ and $j=2,...,q$. Thus, $f(t\tilde{e}_1)=f_{e_1}(t\tilde{e}_1)e_1+f_{\mathcal{Z}}(t\tilde{e}_1) \in \varphi(t)e_1 + \mathcal{Z}$.
\end{proof}
We prove the item 1 of Lemma \ref{lemma_from_arbitrary_to_admissible}.
\begin{sublemma}\label{sublemma2_for_lemma_from_arbitrary_to_admissible}
    The map $\varphi$ is an injective field homomorphism of $\mathbb{F}$, and $f(t\tilde{e}_j)\in \varphi(t)e_j + \mathcal{Z}$.
\end{sublemma}

\begin{proof} First, observe that $\varphi(1)=1$ as $f(\tilde{e}_1)=e_1$. We compute $f([t\tilde{e}_1, s\tilde{e}_j])$ in two ways for $j=2,...,q$ and $t,s\in \mathbb{F}$. First, $f([t\tilde{e}_1, s\tilde{e}_j])=[f(t\tilde{e}_1), f(s\tilde{e}_j)]=\varphi(t)[e_1, f(s\tilde{e}_j)]=\sum_{i=2}^{q} \varphi(t) f_{e_i}(s\tilde{e}_j)[e_1, e_i]$ by Lemma \ref{sublemma1_for_lemma_from_arbitrary_to_admissible}
    
    On the other hand, $f([t\tilde{e}_1, s\tilde{e}_j])=f([ts\tilde{e}_1, \tilde{e}_j])=[f(ts\tilde{e}_1), f(\tilde{e}_j)]=\varphi(ts)[e_1, e_j]$ again by Lemma \ref{sublemma1_for_lemma_from_arbitrary_to_admissible}. By assumption, $[e_1, e_2],...,[e_1,e_q]$ is linearly independent. We see that $\varphi(ts)=\varphi(t)f_{e_j}(s\tilde{e}_j)$ and $\varphi(t)f_{e_i}(s\tilde{e}_j)=0$ for $i\neq j$. The second equation gives $f_{e_i}(s\tilde{e}_j)=0$ (by setting $t=1$), so $f(s\tilde{e}_j)\in f_{e_j}(s\tilde{e}_j)e_j+\mathcal{Z}$. 

     By setting $t=1$ in $\varphi(ts)=\varphi(t)f_{e_j}(s\tilde{e}_j)$, we see that $f_{e_j}(s\tilde{e}_j)=\varphi(s)$, and so $\varphi(ts)=\varphi(t)\varphi(s)$. Since $f$ is additive, $\varphi$ is also additive. This shows that $\varphi$ is a field homomorphism with $\varphi(1)=1$, that is, an injective field homomorphism. Moreover, $f(s\tilde{e}_j)\in\varphi(s)e_j+\mathcal{Z}$ for $j\geq 2$. The case $j=1$ follows from Lemma \ref{sublemma1_for_lemma_from_arbitrary_to_admissible}.
\end{proof}

The key point is that this allows us to show linear independence in $\tilde{\mathcal{V}}$.

\begin{sublemma}\label{lemma_concatenation_linear_independence}
    \begin{enumerate}
        \item $[\tilde{e}_1, \tilde{e}_j]$, $j=2,...,q$, is linearly independent.
        \item $\tilde{e}_1, ..., \tilde{e}_{q}$ is linearly independent, forming an $\mathbb{F}$-basis of $\tilde{\mathcal{V}}$.
    \end{enumerate}
\end{sublemma}

\begin{proof}
    We show the item 1 first. Let $\sum_{j=2}^{q} a_j [\tilde{e}_1, \tilde{e}_j]=0$. Then, $0=f(\sum_{j=2}^{q} a_j [\tilde{e}_1, \tilde{e}_j])=\sum_{j=2}^{q} f([a_j\tilde{e}_1, \tilde{e}_j])=\sum_{j=2}^{q} [f(a_j\tilde{e}_1), f(\tilde{e}_j)]=\sum_{j=2}^{q} \varphi(a_j)[e_1, e_j]$ by Lemma \ref{sublemma1_for_lemma_from_arbitrary_to_admissible}. Since $[e_1, e_j], j=2,...,q,$ is linearly independent, we have $\varphi(a_j)=0$. Since $\varphi$ is injective by Lemma \ref{sublemma2_for_lemma_from_arbitrary_to_admissible}, we have $a_j=0$. This proves the item 1.

    To show the item 2, let $\sum_{j=1}^{q} a_j \tilde{e}_j=0$. Then, $0=[\tilde{e}_1, \sum_{j=1}^{q}a_j \tilde{e}_j]=\sum_{j=2}^{q} a_j[\tilde{e}_1, \tilde{e}_j]$. By the item 1, we have $a_j=0$ for $j=2,...,q$. Now, $a_1 \tilde{e}_1=0$. Since $f(\tilde{e}_1)=e_1\neq 0$, we have $\tilde{e}_1\neq 0$, so $a_1=0$. This proves the item 2.
\end{proof}

Thus, the item 2 of Lemma \ref{lemma_from_arbitrary_to_admissible} follows. Next, we prove the item 3 of Lemma \ref{lemma_from_arbitrary_to_admissible}. 

\begin{sublemma}\label{sublemma3_for_lemma_from_arbitrary_to_admissible}
   For any $x\in \tilde{\mathcal{V}}$ and $t\in \mathbb{F}$, $f(tx)\in \varphi(t)f(x)+\mathcal{Z}$.
\end{sublemma}

\begin{proof}
    Let $t\in \mathbb{F}$ and $x\in \tilde{\mathcal{V}}$ be arbitrary. Then, one can write $x=\sum_{j=1}^{q_1}a_j \tilde{e}_j$, $a_j\in \mathbb{F}$. We have $f(tx)=\sum_{j=1}^{q} f(ta_j\tilde{e}_j)\in \sum_{j=1}^{q} \varphi(ta_j)e_j + \mathcal{Z}=\sum_{j=1}^{q} \varphi(t)\varphi(a_j)e_j+\mathcal{Z}=\varphi(t)f(\sum_{j=1}^{q}a_j \tilde{e}_j)+\mathcal{Z}=\varphi(t)f(x)+\mathcal{Z}$. 
\end{proof}
\begin{proof} [Proof of Lemma \ref{lemma_from_arbitrary_to_admissible}]

The item 1 follows from Lemma \ref{sublemma1_for_lemma_from_arbitrary_to_admissible}. The item 2 follows from Lemma \ref{sublemma2_for_lemma_from_arbitrary_to_admissible}. The item 3 follows from Lemma \ref{sublemma3_for_lemma_from_arbitrary_to_admissible}.
\end{proof}

\subsection{Concatenation} \label{subsection_concatenation}
The $(2n+1)$-dimensional Heisenberg groups and the $(4n+3)$-dimensional quaternionic Heisenberg groups arise as a construction called concatenation (see \cite[Section 3]{Jablonski_moduli}). 

Let $\mathbb{F}$ be any field (including characteristic 2), and let $n\geq 1$ be an integer. 

\begin{definition} \label{definition_concatenation}
    Let $\mathcal{N}$ be a 2-step nilpotent Lie algebra over $\mathbb{F}$.  Let $\mathcal{Z}=\mathfrak{z}(\mathcal{N})$ be the center.  We call $\mathcal{N}$ \textit{a concatenation of admissible 2-step nilpotent Lie algebras} if there are subspaces $\mathcal{V}_m$ of $\mathcal{N}$ such that $\mathcal{N}=\oplus_{m=1}^n \mathcal{V}_m \oplus \mathcal{Z}$ as a vector space over $\mathbb{F}$, and 
\begin{enumerate}
    \item Let $l,m \in \{1,...,n\}$. If $l\neq m$, then $[\mathcal{V}_l,\mathcal{V}_m]=0$.
    \item Each $\mathcal{V}_m$ has an $\mathbb{F}$-basis $e_1^{(m)}, ..., e_{q_m}^{(m)}$ such that $[e_1^{(m)}, e_j^{(m)}]$, $j=2,...,q_m$, is linearly independent.
\end{enumerate}
\end{definition}

\begin{remark}
The subalgebra $\mathcal{V}_m\oplus [\mathcal{V}_m, \mathcal{V}_m]$ is an admissible 2-step nilpotent Lie algebra (see Definition \ref{definition_admissible_Lie_algebras}). Note, a direct sum of admissible 2-step nilpotent Lie algebras is also a concatenation in our definition, unlike the definition in \cite[Section 3]{Jablonski_moduli}.
\end{remark}

\begin{example} The $(2n+1)$-dimensional Heisenberg algebras and $(4n+3)$-dimensional quaternionic Heisenberg algebras are concatenations of admissible 2-step nilpotent Lie algebras (Section \ref{section_definitions}). Any direct sum of admissible 2-step nilpotent Lie algebras is a concatenation of admissible 2-step nilpotent Lie algebras. Indeed, if $\mathcal{N}_m=\mathcal{V}_m\oplus \mathcal{Z}_m$ is an admissible 2-step nilpotent Lie algebra with center $\mathcal{Z}_m$, then $\mathcal{N}:=\oplus_{m=1}^n \mathcal{N}_m=\oplus_{m=1}^n \mathcal{V}_m \oplus \oplus_{m=1}^n \mathcal{Z}_m$ satisfies the conditions in Definition \ref{definition_concatenation}, where the center of $\mathcal{N}$ is $\oplus_{m=1}^n \mathcal{Z}_m$.
\end{example}

%\item $\cap_{m=1}^n [\mathcal{V}_m, \mathcal{V}_m] \neq 0$.

\begin{setting} \label{setting_from_arbitrary_to_concatenation}
    Let $\mathbb{F}$ be a field. Let $\mathcal{N}$ be a 2-step nilpotent Lie algebra over $\mathbb{F}$. Let $\mathcal{Z}=\mathfrak{z}(\mathcal{N})$ be the center of $\mathcal{N}$. Suppose that $\mathcal{N}$ is a concatenation of admissible 2-step nilpotent Lie algebras, so that $\mathcal{N}=\oplus_{m=1}^n \mathcal{V}_m \oplus \mathcal{Z}$ as in Definition \ref{definition_concatenation}. Let $\mathcal{V}=\oplus_{m=1}^n \mathcal{V}_m$. Let $\tilde{\mathcal{N}}$ be an arbitrary 2-step nilpotent Lie algebra over $\mathbb{F}$. Let $f:\tilde{\mathcal{N}}\to \mathcal{N}$ be a Lie ring \textit{isomorphism}.
\end{setting}

\begin{remark}
    This is an analogue of Setting \ref{setting_technical_lemma}. This time, we assume $f$ is a Lie ring isomorphism.
\end{remark}

\begin{notation} \label{notation__from_arbitrary_to_concatenation}
Under Setting \ref{setting_from_arbitrary_to_concatenation}, let $\mathcal{N}=\mathcal{V}\oplus \mathcal{Z}=\oplus_{m=1}^n \mathcal{V}_m \oplus \mathcal{Z}$ with notations as in Definition \ref{definition_concatenation}. Let $\tilde{e}_j^{(m)}=f^{-1}(e_j^{(m)})\in \tilde{\mathcal{N}}$. Let us enumerate $e_1^{(1)},...,e_{q_1}^{(1)}, e_1^{(2)},..., e_{q_2}^{(2)},..., e_1^{(n)}, ..., e_{q_n}^{(n)}$ as $e_1, ..., e_{q_1}, e_{q_1+1}, ..., e_{q_1+q_2}, ..., e_{Q}$, where $Q=q_1+...+q_n$. Write $\tilde{e}_j=f^{-1}(e_j)$. Let $\tilde{\mathcal{V}}_m=\mathbb{F}\textrm{-span}\{\tilde{e}_1^{(m)},...,\tilde{e}_{q_m}^{(m)}\}$. Let $\tilde{\mathcal{V}}=\sum_{m=1}^n \tilde{\mathcal{V}}_m$. Given $X\in \tilde{\mathcal{N}}$, write $f(X)=\sum_{j=1}^Q f_{e_j}(X)e_j+f_{\mathcal{Z}}(X)$, where $f_{\mathcal{Z}}(X)\in \mathcal{Z}$. 
\end{notation}
Observe, any element $x\in \tilde{\mathcal{V}}$ has the form $x=\sum_{j=1}^Q a_j\tilde{e}_j$ for some $a_j\in \mathbb{F}$ by definition.
\begin{lemma} \label{lemma_from_arbitrary_to_concatenation}
    Under Setting \ref{setting_from_arbitrary_to_concatenation} with Notation \ref{notation__from_arbitrary_to_concatenation}, let $\varphi_m(t)=f_{e_1^{(m)}}(t\tilde{e}_1^{(m)})$ for each $m\in \{1,...,n\}$. Then, for each $m\in \{1,...,n\}$, we have 
    \begin{enumerate}
        \item Let $l,m \in \{1,...,n\}$. If $l\neq m$, then $[\tilde{\mathcal{V}}_l, \tilde{\mathcal{V}}_m]=0$.
        \item $[\tilde{e}_1^{(m)}, \tilde{e}_j^{(m)}], j\in \{2,...,q_m\}$ is linearly independent, and $\tilde{e}_1^{(m)},...,\tilde{e}_{q_m}^{(m)}$ is linearly independent.
        \item $\varphi_m$ is an injective field homomorphism of $\mathbb{F}$, and $f(t\tilde{e}_j^{(m)})\in \varphi_m(t)e_j^{(m)}+\mathcal{Z}$ for $j\in \{1,...,q_m\}$.
        \item $f(\tilde{\mathcal{V}}_m) \subseteq \mathcal{V}_m\oplus\mathcal{Z}$ for any $m=1,...,n$. For any $x\in \tilde{\mathcal{V}}_m$ and $t\in \mathbb{F}$, $f(tx)\in \varphi_m(t)f(x)+\mathcal{Z}$.
        \item If $\mathbb{F}=\mathbb{R}$, then $f(\sum_{j=1}^Q a_j \tilde{e}_j)\in \sum_{j=1}^Q a_j e_j +\mathcal{Z}$ for any $\sum_{j=1}^Q a_j \tilde{e}_j\in \tilde{\mathcal{V}}$. For any $t\in \mathbb{R}$ and $x\in \tilde{\mathcal{V}}$, $f(tx)\in tf(x)+\mathcal{Z}$.
    \end{enumerate}
\end{lemma}
\begin{remark}
    This is an analogue of Lemma \ref{lemma_from_arbitrary_to_admissible}. The field homomorphism may not be surjective, as Example \ref{example_arising_from_field_of_rational_functions} shows.
\end{remark}

By the following sublemma, we can reduce the proof of Lemma \ref{lemma_from_arbitrary_to_concatenation} to Lemma \ref{lemma_from_arbitrary_to_admissible}. 

\begin{sublemma} Under Setting \ref{setting_from_arbitrary_to_concatenation} with Notation \ref{notation__from_arbitrary_to_concatenation}, we have
\label{properties_of_tilde}
    \begin{enumerate} 
        \item Let $l,m \in \{1,...,n\}$. If $l\neq m$, then $[\tilde{\mathcal{V}}_l, \tilde{\mathcal{V}}_m]=0$.
        \item $\cap_{j=q_1+1}^Q C(e_j)=\mathcal{V}_1 \oplus \mathcal{Z}$, where $C(e_j)$ is the centralizer of $e_j$ in $\mathcal{N}$.
        \item $f(\tilde{\mathcal{V}}_m) \subseteq \mathcal{V}_m\oplus\mathcal{Z}$ for any $m=1,...,n$.
    \end{enumerate}
\end{sublemma}

\begin{proof}
    We prove the item 1. If $l\neq m$, then $f([\tilde{e}_i^{(l)}, \tilde{e}_j^{(m)}])=[e_i^{(l)},e_j^{(m)}]=0$ by Setting \ref{setting_from_arbitrary_to_concatenation}. Since $f$ is injective, we have $[\tilde{e}_i^{(l)}, \tilde{e}_j^{(m)}]=0$. 
    
    Next, we prove the item 2. Since $[\mathcal{V}_1, \mathcal{V}_m]=0$ for $m\geq 2$, we have one inclusion. To show the other inclusion, let $x\in \cap_{j=q_1+1}^Q C(e_j)$, and write $x=\sum_{m=1}^n \sum_{j=1}^{q_m} a_j^{(m)} e_j^{(m)}+z$, where $z\in \mathcal{Z}$. Then, for any $m\geq 2$, $0=[e_1^{(m)}, x]=\sum_{j=2}^{q_m}a_j^{(m)} [e_1^{(m)}, e_j^{(m)}].$ By linear independence (the item 2 in Definition \ref{definition_concatenation}), we have $a_j^{(m)}=0$ for $j=2,...,q_m$ and for any $m\geq 2$. Now, $x=\sum_{m=2}^n a_1 e_1^{(m)} + \sum_{j=1}^{q_1} a_j^{(1)} e_j^{(1)}+z$. Then, for any $m\geq 2$, $0=[x, e_2^{(m)}]=a_1^{(m)}[e_1^{(m)},e_2^{(m)}]$. By the item 2 in Definition \ref{definition_concatenation}, $[e_1^{(m)},e_2^{(m)}]\neq 0$. Thus, $a_1^{(m)}=0$. This shows $x\in \mathcal{V}_1\oplus \mathcal{Z}$. 

    Finally, we show the item 3. Without loss of generality, we may assume $m=1$. Let $t\in \mathbb{F}$ and $i=1,..., q_1$. By the item 1, for any $j\geq q_1+1$, we have $[\tilde{e}_i, \tilde{e}_j]=0$. Thus, $[t\tilde{e}_i, \tilde{e}_j]=0$. Hence, $0=f([t\tilde{e}_i, \tilde{e}_j])=[f(t\tilde{e}_i), f(\tilde{e}_j)]=[f(t\tilde{e}_i), e_j]$. We have $f(t\tilde{e}_i)\in \cap_{j=q_1+1}^Q C(e_j)=\mathcal{V}_1\oplus \mathcal{Z}$ by the item 2. Since $f$ is additive and $\tilde{e}_1,...,\tilde{e}_{q_1}$ span $\tilde{\mathcal{V}}_1$, this shows $f(\tilde{\mathcal{V}}_1)\subseteq \mathcal{V}_1\oplus \mathcal{Z}.$
\end{proof}

We are ready to prove Lemma \ref{lemma_from_arbitrary_to_concatenation}.

\begin{proof} [Proof of Lemma \ref{lemma_from_arbitrary_to_concatenation}.] The item 1 follows from Lemma \ref{properties_of_tilde}. Let $\mathcal{N}_m=\mathcal{V}_m\oplus \mathcal{Z}$.
    Let $\tilde{\mathcal{N}}_m=\tilde{\mathcal{V}_m} + \tilde{\mathcal{Z}}$, where $\tilde{\mathcal{Z}}$ is the center of $\tilde{\mathcal{N}}$. Then, $\mathcal{N}_m, \tilde{\mathcal{N}_m}$ are 2-step nilpotent Lie algebras over $\mathbb{F}$. By Lemma \ref{center_goes_to_center}, $f(\tilde{\mathcal{Z}})=\mathcal{Z}$. By Lemma \ref{properties_of_tilde}, we have $f(\tilde{\mathcal{V}}_m)\subseteq \mathcal{V}_m\oplus \mathcal{Z}$. Thus, $f|_{\tilde{\mathcal{N}}_m}:{\tilde{\mathcal{N}}_m}\to \mathcal{N}_m$ is a Lie ring homomorphism such that $f(\tilde{e}_j^{(m)})=e_j^{(m)}$. 

    To apply Lemma \ref{lemma_from_arbitrary_to_admissible}, we show $\mathcal{N}_m$ is admissible relative to $e_j^{(m)}, j\in \{1,...,q_m\}$. By Setting \ref{setting_from_arbitrary_to_concatenation}, the only condition we must check is that $\mathcal{Z}=\mathfrak{z}(\mathcal{N}_m)$. As $\mathcal{Z}$ is the center of $\mathcal{N}$, it is contained in $\mathfrak{z}(\mathcal{N}_m)$. Conversely, if $z\in \mathfrak{z}(\mathcal{N}_m)$, then write it as $z= \sum_{j=1}^{q_m} a_j e_j^{(m)}+z'$, where $z'\in \mathcal{Z}$. Then, $0=[e_1^{(m)},z]=\sum_{j=2}^{q_m} a_j [e_1^{(m)}, e_j^{(m)}]$. By the item 2 of Definition \ref{definition_concatenation}, $a_j^{(m)}=0$ for $j\in \{2,...,q_m\}$. Then, $0=[z,e_2^{(m)}]=a_1[e_1^{(m)},e_2^{(m)}]$ and again by linear independence, we have $a_1^{(m)}=0$. This shows $\mathcal{Z}=\mathfrak{z}(\mathcal{N}_m)$.

    Thus, Lemma \ref{lemma_from_arbitrary_to_admissible} applies to $f|_{\tilde{\mathcal{N}}_m}:\tilde{\mathcal{N}}_m \to \mathcal{N}_m$, and the items 2, 3, and 4 follow. If $\mathbb{F}=\mathbb{R}$, then the only nontrivial field homomorphism of $\mathbb{R}$ is the identity (see \cite[Theorem 14.4.1]{Kuczma}). Thus, the first part of the item 5 follows from the item 3. Any $x\in \tilde{\mathcal{V}}$ has the form $x=\sum_{j=1}^Q a_j\tilde{e}_j$ for some $a_j\in \mathbb{R}$. For any $t\in \mathbb{R}$, $f(tx)\in \sum_{j=1}^Q ta_j \tilde{e}_j +\mathcal{Z}=tf(x)+\mathcal{Z}$. This completes the proof.
\end{proof}

The following theorem is our main technical theorem.

\begin{theorem}\label{theorem_concatenation}
    Let $\mathbb{F}$ be any field. Let $\mathcal{N}$ be a concatenation of admissible 2-step nilpotent Lie algebras as in Definition \ref{definition_concatenation}. We use the notations in Definition \ref{definition_concatenation}. Assume $\cap_{m=1}^n [\mathcal{V}_m, \mathcal{V}_m] \neq 0$. Let $\mathcal{V}=\oplus_{m=1}^n \mathcal{V}_m$.
    
    Then,
    
    \begin{enumerate}
        \item Any Lie ring automorphism of $\mathcal{N}$ is a product of a central automorphism and $\varphi$-linear Lie ring automorphism of $\mathcal{N}$ for some field automorphism $\varphi$ of $\mathbb{F}$.
        \item Assume that the center $\mathcal{Z}=\mathfrak{z}(\mathcal{N})$ coincides with $[\mathcal{N},\mathcal{N}]$. Let $f$ be a Lie ring automorphism of $\mathcal{N}$ such that $f(\mathcal{V})=\mathcal{V}$. Then, there is a field automorphism $\varphi$ of $\mathbb{F}$ such that $f$ is $\varphi$-linear.
    \end{enumerate}

    Moreover, if $\mathbb{F}=\mathbb{R}$, the same result holds even when $\cap_{m=1}^n [\mathcal{V}_m, \mathcal{V}_m] = 0$.
\end{theorem}

\begin{remark} \label{remark_why_nonempty_intersection}
    The condition $\cap_{m=1}^n [\mathcal{V}_m, \mathcal{V}_m] \neq 0$ is only used to make sure $\varphi_m=\varphi_1$ for each $m$ in Lemma \ref{lemma_from_arbitrary_to_concatenation}. Without this condition, $\mathcal{N}$ can be decomposable.
\end{remark}

\begin{proof}  
    We show the item 2 first. Setting \ref{setting_from_arbitrary_to_concatenation} is satisfied with $\tilde{\mathcal{N}}=\mathcal{N}$. Use the same notations as in Notation \ref{notation__from_arbitrary_to_concatenation}. Then, $f(\mathcal{V})=\mathcal{V}$ implies that $\tilde{\mathcal{V}}_m\subseteq \mathcal{V}$, $f(\tilde{\mathcal{V}}_m)\subseteq \mathcal{V}_m$, and $f(tx)=\varphi_m(t)f(x)$ for any $t\in \mathbb{F}$ and $x\in \tilde{\mathcal{V}}_m$ by Lemma \ref{lemma_from_arbitrary_to_concatenation}.

    We show $\mathcal{V}=\oplus_{m=1}^n \tilde{\mathcal{V}}_m$. First, we show $\oplus_{m=1}^n \tilde{\mathcal{V}}_m=\sum_{m=1}^n \tilde{\mathcal{V}}_m$. Let $\sum_{m=1}^n x_m=0$, where $x_m\in \tilde{\mathcal{V}}_m$. Then, $\sum_{m=1}^n f(x_m)=0$, $f(x_m)\in \mathcal{V}_m$. Thus, $f(x_m)=0$. Since $f$ is injective, $x_m=0$. Thus, $\oplus_{m=1}^n \tilde{\mathcal{V}}_m$ is a subspace of dimension $Q=q_1+...+q_m$ of $\mathcal{V}$. Since $\dim \mathcal{V}=Q$, we have $\mathcal{V}=\oplus_{m=1}^n \tilde{\mathcal{V}}_m$. 

    Next, we show $\varphi_m$ is a field automorphism of $\mathbb{F}$. As $\mathcal{N}=\mathcal{V}\oplus\mathcal{Z}=\oplus_{m=1}^n \tilde{\mathcal{V}}_m\oplus \mathcal{Z}$, $\mathcal{N}$ is a concatenation of admissible 2-step nilpotent Lie algebras with respect to $\oplus_{m=1}^n \tilde{\mathcal{V}}_m$ by Lemma \ref{lemma_from_arbitrary_to_concatenation}. Thus, Lemma \ref{lemma_from_arbitrary_to_concatenation} applies to $f^{-1}$ as well. Thus, for each $m\in \{1,...,n\}$, there is a field homomorphism $\psi_m$ of $\mathbb{F}$ such that $f^{-1}(te_j^{(m)})=\psi_m(t)e_j^{(m)}$. Then, $te_1^{(m)}=f(f^{-1}(te_1^{(m)}))=\varphi_m(\psi_m(t))e_1^{(m)}$. Thus, $\varphi_m$ is a field automorphism of $\mathbb{F}$. 
    
    Let $\varphi=\varphi_1$. We show $\varphi=\varphi_m$ for any $m$. This is the part we use $\cap_{m=1}^n [\mathcal{V}_m, \mathcal{V}_m]\neq 0$. Let $0\neq z\in \cap_{m=1}^n [\mathcal{V}_m, \mathcal{V}_m]$. Fix $m=1,...,n$, and write $z=\sum_{j=1}^l [x_j,y_j]$, where $x_j,y_j\in \mathcal{V}_m$. Since $f^{-1}$ is $\psi_m$-linear on each $\mathcal{V}_m$, we have $f^{-1}(tz)=f^{-1}(\sum_{j=1}^l [tx_j,y_j])=\sum_{j=1}^l [f^{-1}(tx_j), f^{-1}(y_j)]=\psi_m(t)\sum_{j=1}^l[f^{-1}(x_j), f^{-1}(y_j)]=\psi_m(t)f^{-1}(z)$. However, $f^{-1}(tz)$ does not depend on $m$. Thus, $\psi_m(t)f^{-1}(z)=\psi_1(t)f^{-1}(z)$. Since $z\neq 0$, we have $f^{-1}(z)\neq 0$, and $\psi_m(t)=\psi_1(t)$. Thus, $\varphi_m=\psi_m^{-1}=\psi_1^{-1}=\varphi_1=\varphi$. 

Therefore, for any $x=\sum_{m=1}^n x_m$, $x_m\in \tilde{\mathcal{V}}_m$, and for any $t\in \mathbb{F}$, we have $f(tx)=\sum_{m=1}^n f(tx_m)=\sum_{m=1}^n \varphi(t)f(x_m)=\varphi(t)f(x)$. This shows $f|_{\mathcal{V}}$ is $\varphi$-linear. 

We show $f|_{\mathcal{Z}}$ is also $\varphi$-linear. Since $[\mathcal{N},\mathcal{N}]=\mathcal{Z}$, any $z\in \mathcal{Z}$ has the form $z=\sum_{j=1}^l [x_j, y_j] \in \mathcal{Z}$, where $x_j, y_j\in \mathcal{V}$. Then, for any $t\in \mathbb{F}$, $f(tz)=\sum_{j=1}^l [f(tx_j), f(y_j)]=\varphi(t)\sum_{j=1}^l[f(x_j),f(y_j)]=\varphi(t)f(z)$. Thus, $f|_{\mathcal{Z}}$ is $\varphi$-linear, so $f$ is $\varphi$-linear.

We show the item (1). By Proposition \ref{abelian_factors_2-step} (2), which is independently shown in Section \ref{section_abelian_factors}, we may assume $\mathcal{Z}=[\mathcal{N}, \mathcal{N}]$, and we can use the item (2). We show any Lie ring automorphism $f$ of $\mathcal{N}$ is a product of a central automorphism and a $\varphi$-linear Lie ring automorphism of $\mathcal{N}$. By Proposition \ref{semidirect_product_for_2-step_nilpotent_Lie_algebras_over_fields}, which is independently shown in Section \ref{section_semidirect_product}, $f=\mu\circ \overline{f}$, where $\mu$ is a central automorphism and $\overline{f}$ is a Lie ring automorphism with $\overline{f}(\mathcal{V})=\mathcal{V}$. Thus, $\overline{f}$ is a $\varphi$-linear Lie ring automorphism for some field automorphism $\varphi$ of $\mathbb{F}$.

If $\mathbb{F}=\mathbb{R}$, then $\varphi_m(t)=t$, as the identity is the only field automorphism of $\mathbb{R}$. Thus, even when $\cap_{m=1}^n [\mathcal{V}_m, \mathcal{V}_m]=0$, we have $\varphi_m=\varphi=\textrm{id}_{\mathbb{R}}$, and the same result holds. This completes the proof.
\end{proof}

\begin{corollary} \label{corollary_concatenation}
    Let $\mathbb{F}=\mathbb{R}$ or $\mathbb{C}$. Let $\mathcal{N}$ be a concatenation of admissible 2-step nilpotent Lie algebras as in Definition \ref{definition_concatenation}. We use the notations in Definition \ref{definition_concatenation}. Let $N$ be a simply connected 2-step nilpotent Lie group with Lie algebra $\mathcal{N}$.
    \begin{enumerate}
        \item Let $\mathbb{F}=\mathbb{R}$. Then, any abstract group automorphism of $N$ is a product of a central automorphism and a Lie group automorphism.
        \item Let $\mathbb{F}=\mathbb{C}$. Assume that $\cap_{m=1}^n [\mathcal{V}_m, \mathcal{V}_m] \neq 0$ (cf. Remark \ref{remark_why_nonempty_intersection}). Then, every abstract group automorphism of $N$ is standard.
    \end{enumerate}
\end{corollary}
\begin{proof}
    This follows from Theorem \ref{theorem_concatenation}, Lemma \ref{reduction_to_the_Lie_algebras}, and the fact that $\varphi$ must be $\textrm{id}_{\mathbb{R}}$ for $\mathbb{F}=\mathbb{R}$.
\end{proof}
Given a 2-step nilpotent Lie algebra $\mathcal{N}$ over a field $\mathbb{F}$ of characteristic not 2, let $N$ be $\mathcal{N}$ as a set with the multiplication $xy:=x+y+\frac{1}{2}[x,y]$. Since $\mathcal{N}$ is 2-step nilpotent, $(xy)z=x(yz)$. The inverse map is given by $x^{-1}:=-x$. The group $N$ is called the \textit{2-step nilpotent affine algebraic group over $\mathbb{F}$ associated to} $\mathcal{N}$.

  If $\mathbb{F}=\mathbb{R}$ or $\mathbb{C}$, then this definition recovers a simply connected 2-step nilpotent Lie group whose Lie algebra is $\mathcal{N}$ by the Baker--Campbell--Hausdorff formula.

If $f$ is a Lie ring automorphism of $\mathcal{N}$, then the same map is clearly an abstract group automorphism of $N$. Conversely, if $F$ is an abstract group automorphism of $N$, then $F[x,y]=[F(x),F(y)]$ because $xyx^{-1}y^{-1}=[x,y]$. Noting that $xz=x+z$ for $z\in [\mathcal{N}, \mathcal{N}]$, we have $F(x+y)=F(xy)F(-\frac{1}{2}[x,y])$. Also, $F(x)+F(y)=F(x)F(y)(-\frac{1}{2}[F(x),F(y)])=F(xy)(-\frac{1}{2}F[x,y])$. Setting $x=y$, one obtains $F(2x)=2F(x)$, which in turn gives $F(\frac{1}{2}x)=\frac{1}{2}F(x)$ as $2$ is invertible. Hence, $F(x+y)=F(x)+F(y)$, i.e., $F$ is a Lie ring automorphism of $\mathcal{N}$. 
\begin{remark} \label{remark_Baer_correspondence}
    Similarly, if $ \mathcal{N}, \mathcal{N}'$ are 2-step nilpotent Lie algebras over $\mathbb{F}$ and $N, N'$ are the groups associated to $\mathcal{N}, \mathcal{N}'$, respectively (see the definition above), then a map from $N\to N'$ is an abstract group isomorphism if and only if the same set map $\mathcal{N}\to \mathcal{N}'$ is a Lie ring isomorphism.
\end{remark}
 This is known as Baer correspondence \cite{Baer}. In this case, the exponential map $\mathcal{N}\to N$ is the identity map, and a field isomorphism of Lie algebras $\overline{\varphi}:\mathcal{N}\to \mathcal{N}^\varphi$ gives rise to the group isomorphism $\hat{\varphi}=\overline{\varphi}:N\to N^\varphi$, which we call \textit{the field isomorphism of groups} (cf. Section \ref{section_introduction_complex_case}). 

\begin{corollary} \label{corollary_field_of_char_not_2} Let $\mathbb{F}$ be any field of characteristic not 2. Let $\mathcal{N}$ be a concatenation of admissible 2-step nilpotent Lie algebras as in Definition \ref{definition_concatenation}. We use the notations in Definition \ref{definition_concatenation}. Assume $\cap_{m=1}^n [\mathcal{V}_m, \mathcal{V}_m] \neq 0$. Let $N$ be the 2-step nilpotent affine algebraic group over $\mathbb{F}$ associated to $\mathcal{N}$, as defined above.

Then, for every abstract group automorphism $F$ of $N$, there exist a central automorphism $\mu$, a field automorphism $\varphi$ of $\mathbb{F}$, and a morphism $\overline{F}:N^\varphi \to N$ such that $F=\mu \circ \overline{F} \circ \hat{\varphi}$, where $\hat{\varphi}:N\to N^\varphi$ is the induced field isomorphism of groups. Here, $N^\varphi$ is the group associated to $\mathcal{N}^\varphi$.
\end{corollary}
\begin{proof}
    By the discussion above (\cite{Baer}), $F$ is a Lie ring automorphism of $\mathcal{N}$. By Theorem \ref{theorem_concatenation}, $F=\mu\circ f$, where $\mu$ is a central automorphism and $f$ is a $\varphi$-linear Lie ring automorphism of $\mathcal{N}$ for some $\varphi\in \textrm{Aut}(\mathbb{F})$. By Proposition \ref{existence_of_semilinear_Lie_ring_automorphism} (2), $f=\overline{f}\circ \overline{\varphi}$, where $\overline{f}:\mathcal{N}^\varphi\to \mathcal{N}$ is a Lie algebra isomorphism and $\overline{\varphi}:\mathcal{N}\to \mathcal{N}^\varphi$ is a field isomorphism. This $\overline{F}:=\overline{f}$ is a morphism, and $\hat{\varphi}:=\overline{\varphi}$ is a field isomorphism of groups (see Remark \ref{remark_Baer_correspondence}).
\end{proof}

\section{Proof of Theorem \ref{generic_condition_introduction}, Theorem \ref{maintheorem_criterion}, and Cases (1) and (2) of Theorem \ref{maintheorem_over_C}} \label{section_generic}

Let $\mathbb{F}=\mathbb{R}$ or $\mathbb{C}$ in this section. We begin with the proof of Theorem \ref{maintheorem_criterion} and Theorem \ref{maintheorem_over_C} for the case of (2).

\begin{proof}[Proof of Theorem \ref{maintheorem_criterion}] This follows from Corollary \ref{corollary_concatenation} with $\mathbb{F}=\mathbb{R}$ and $n=1$.
\end{proof}

\begin{proof}[Proof of the case (2) of Theorem \ref{maintheorem_over_C}]
This follows from Corollary \ref{corollary_concatenation} with $\mathbb{F}=\mathbb{C}$ and $n=1$. 
\end{proof}

Next, we will show that the condition of $\mathcal{N}$ being admissible is either generic or empty. To make it more precise, we introduce the space of 2-step nilpotent Lie algebras. See \cite[Section 2, 3]{Console_Sergio_Fino_Evangelia} for more details. Any simply connected 2-step nilpotent Lie group $N$ is uniquely determined by its Lie algebra $\mathcal{N}$ up to isomorphism. We say that a 2-step nilpotent Lie algebra $\mathcal{N}$ is \textit{of type $(p,q)$} if $\dim_{\mathbb{F}} [\mathcal{N}, \mathcal{N}]=p$ and $\dim_{\mathbb{F}} \mathcal{N}=p+q$. To consider the space of 2-step nilpotent Lie algebras of type $(p,q)$, we fix the vector space $\mathbb{F}^q\oplus \mathbb{F}^p$. By choosing a basis, we may identify $\mathcal{N}$ as $(\mathbb{F}^q\oplus \mathbb{F}^p, [\cdot,\cdot])$ with $\mathbb{F}^p$ the commutator ideal, where $[\cdot,\cdot]$ is the corresponding Lie bracket. Since the commutator ideal lies in the center, the non-zero part of the Lie bracket is given by a surjective linear map $\Lambda^2(\mathbb{F}^q) \to \mathbb{F}^p$. Let $\mathcal{V}(p,q)=\Lambda^2(\mathbb{F}^q)^*\otimes \mathbb{F}^p$. Conversely, any linear map $\rho\in \mathcal{V}(p,q)$ determines a 2-step nilpotent Lie algebra $\mathcal{N}_{\rho}=(\mathbb{F}^q\oplus \mathbb{F}^p, [\cdot, \cdot]_{\rho})$, and it is of type $(p,q)$ if $\rho$ is surjective. The set of 2-step nilpotent Lie algebras of type $(p,q)$ forms a Zariski open, dense subset $V_{p,q}^o$ of $\mathcal{V}(p,q)$.  

Given $q\geq 2$, we must have $1\leq p\leq \frac{1}{2}q(q-1)$ (\cite[p.46]{Eberlein_moduli_2-step}). For each $\rho\in V_{p,q}^o$, let $N_\rho$ be a simply connected 2-step nilpotent Lie group whose Lie algebra is $\mathcal{N}_\rho$. 

\begin{theorem} \label{generic_condition_precise}
    Let $p\geq 1, q\geq 2$ be integers such that $q-1\leq p\leq \frac{1}{2}q(q-1)$. Let $\mathcal{S}$ be the set of $\rho\in V_{p,q}^o$ such that $\ker ad_{\rho} e_1=\mathbb{F}\textrm{-span}\{e_1\}\oplus \mathbb{F}^p$, where $e_1, ..., e_q$ is the standard basis of $\mathbb{F}^q$. Then,
    \begin{enumerate}
        \item For any $\rho\in \mathcal{S}$,  the center coincides with $\mathbb{F}^p$.
        \item There exists a Zariski open, dense subset $\mathcal{O}$ of $\mathcal{V}(p,q)$ such that $\mathcal{O}\cap V_{p,q}^o \subseteq \mathcal{S}$.
        \item Let $\mathbb{F}=\mathbb{R}$. Then for any $\rho\in \mathcal{O}\cap V_{p,q}^o$, $N_\rho$ satisfies partial automatic continuity.
        \item Let $\mathbb{F}=\mathbb{C}$. Then for any $\rho\in \mathcal{O}\cap V_{p,q}^o$, any abstract group automorphism of $N_\rho$ is standard.
    \end{enumerate}
    
\end{theorem}

\begin{proof}
    Given any $\rho\in \mathcal{V}(p,q)$, let $ad_{\rho}x(y)=[x, y]_{\rho}$, $x, y\in \mathbb{F}^q\oplus \mathbb{F}^p$. We first show if $\rho\in \mathcal{S}$, then the center coincides with $\mathbb{F}^p$. Suppose $Z$ is in the center. Then $Z\in \ker ad_{\rho} e_1=\mathbb{R}\textrm{-span}\{e_1\}\oplus \mathbb{F}^p$, so $Z=a_1e_1+v$, $v\in \mathbb{F}^p$. Since $e_1$ cannot be in the center, $a_1=0$, so $Z\in \mathbb{F}^p$. 

Next, we show $\mathcal{S}$ contains a Zariski open set. Note, $ad_{\rho} e_1|_{\mathbb{F}^q}\in (\mathbb{F}^q)^*\otimes \mathbb{F}^p$, and $\sigma:\mathcal{V}(p,q)\to (\mathbb{F}^q)^*\otimes \mathbb{F}^p$, $\sigma(\rho)=ad_{\rho}e_1|_{\mathbb{F}^q}$, is a polynomial map. Let $\{\overline{e}_1,...,\overline{e}_p\}$ be the standard basis of $\mathbb{F}^p$. Let $\det:(\mathbb{F}^q)^*\otimes \mathbb{F}^p\to \mathbb{F}$ be the map that assigns the determinant of the $(q-1)\times (q-1)$ submatrix corresponding to $e_2,...,e_q$ and $\overline{e}_1,...,\overline{e}_{q-1}$. This is well defined since $p\geq q-1$ and also a polynomial map. Let $\mathcal{O}=(\det\circ \sigma)^{-1}(\mathbb{F}\setminus\{0\})$. If $\rho\in \mathcal{O}\cap V_{p,q}^o$, then $\dim \ker ad_{\rho} e_1 \leq p+1$. As $\ker ad_{\rho} e_1$ always contains $\mathbb{F}\textrm{-span}\{e_1\}\oplus \mathbb{F}^p$, this implies $\rho \in \mathcal{S}$. 

We show that $\mathcal{O}\cap V_{p,q}^o$ is non-empty. Let $\mathbb{F}^q \oplus \Lambda^2 \mathbb{F}^q$ be a free 2-step nilpotent Lie algebra. The commutator ideal has dimension $D=\frac{1}{2}q(q-1)$. As $p\geq q-1$, choose $p$ many relations including $[e_1, e_2],...,[e_1,e_q]$, and let $\mathcal{N}$ be the quotient by the rest of the relations. There is $\rho\in \mathcal{O}\cap V_{p,q}^o$ that is isomorphic to $\mathcal{N}$. The item 2 follows.

If $\rho\in \mathcal{S}$, since $\mathbb{F}^{p+q}=\mathbb{F}\textrm{-span}\{e_2, ..., e_q\}\oplus \ker ad_{\rho} e_1$, it follows that $[e_1, e_2], ..., [e_1, e_q]$ is linearly independent. Thus, by Corollary \ref{corollary_concatenation} with $n=1$, the items 3 and 4 follow.
\end{proof}

\begin{proof}[Proof of Theorem \ref{generic_condition_introduction} and the case (1) of Theorem \ref{maintheorem_over_C}]
    By Remark \ref{remark_reason_for_bounds}, the condition $q-1\leq p$ is the same as $\dim[N,N]\geq \frac{\dim N -1}{2}$. Thus, Theorem \ref{generic_condition_precise} (3), (4) imply the results. 
\end{proof}

\section{Proof of Theorem \ref{pac_up_to_dim_6_introduction}}\label{section_counterexamples}

Write $\textrm{Heis}_{2n+1}$ for the $(2n+1)$-dimensional Heisenberg group and $\textrm{heis}_{2n+1}$ for the $(2n+1)$-dimensional Heisenberg algebra.  Let $\mathbb{R}[\epsilon]:=\mathbb{R}[X]/\langle X^2 \rangle$ be the ring of dual numbers, where $\epsilon$ denotes $X$ in the quotient. Then, $\mathbb{R}[\epsilon]=\mathbb{R}\textrm{-span}\{1,\epsilon\}$ is a 2-dimensional $\mathbb{R}$-algebra. 

\begin{example}\label{def_of_counter_example}
Consider $\textrm{Heis}_3(\mathbb{R}[\epsilon])$ with Lie algebra $\textrm{heis}_3(\mathbb{R}[\epsilon])$. Concretely, one can realize $\mathbb{R}[\epsilon]$ as a subalgebra of $M_2(\mathbb{R})$ and $\textrm{heis}_3(\mathbb{R}[\epsilon])\subseteq M_3(M_2(\mathbb{R}))\cong M_6(\mathbb{R})$. Let $X, Y, Z$ denote the usual basis of $\textrm{heis}_3(\mathbb{R})$ with $[X,Y]=Z$. Then, $\textrm{heis}_3(\mathbb{R}[\epsilon])=\{aX+bY+cZ\;|\;a,b,c\in \mathbb{R}[\epsilon]\}$ has a $\mathbb{R}$-basis $X_1=X, X_2=Y, X_3=\epsilon Y, X_4=-\epsilon X, Z_1=Z, Z_2=\epsilon Z$. Then, $[X_1, X_2]=Z_1, [X_1, X_3]=Z_2, [X_1, X_4]=0, [X_2, X_3]=0, [X_2, X_4]=Z_2, [X_3, X_4]=0$. 
\end{example}

For any nilpotent Lie algebra $\mathcal{N}$, there is a unique simply connected nilpotent Lie group whose Lie algebra is $\mathcal{N}$ up to isomorphism. Thus, we can state the classification at the Lie algebra level, as shown in Table \ref{table_low_dimensions}. Let $\mathcal{N}_5:=\mathbb{R}\text{-span}\{X_1, X_2, X_3, Z_1, Z_2\}$ with $[X_1, X_2]=Z_1, [X_1, X_3]=Z_2$. Let $\mathcal{N}_6:=\mathbb{R}\textrm{-span}\{X_1, X_2, X_3\}\oplus \mathbb{R}\textrm{-span}\{Z_1, Z_2, Z_3\}$ with $[X_1, X_2]=Z_3, [X_2,X_3]=Z_1, [X_3,X_1]=Z_2$.

\begin{table}[ht]
\begin{center}
\caption{Classification of 2-step nilpotent Lie algebras of dimensions 6 or less}
\label{table_low_dimensions}
\begin{tabular}{ |c|c| } 
 \hline
 $\dim 3, 4, 5$ & $\textrm{heis}_3(\mathbb{R})$, $\textrm{heis}_3(\mathbb{R})\oplus \mathbb{R}$, $\textrm{heis}_3(\mathbb{R})\oplus \mathbb{R}^2$, $\textrm{heis}_5(\mathbb{R})$, $\mathcal{N}_5$ \\ 
 \hline
 $\dim 6$ with an abelian factor & $\textrm{heis}_3(\mathbb{R})\oplus \mathbb{R}^3$, $\textrm{heis}_5(\mathbb{R})\oplus \mathbb{R}$, $\mathcal{N}_5\oplus \mathbb{R}$,  \\ 
 \hline
 $\dim 6$ without an abelian factor & $\textrm{heis}_3(\mathbb{R})\oplus \textrm{heis}_3(\mathbb{R})$, $\textrm{heis}_3(\mathbb{C})$, $\textrm{heis}_3(\mathbb{R}[\epsilon])$, $\mathcal{N}_6$ \\ 
 \hline
\end{tabular}
\end{center}
\end{table}

Let us explain the table briefly. For dimensions 3,4, and 5, one can consult \cite[p.22]{Console_Sergio_Fino_Evangelia} (see also \cite[Theorem 10]{Homolya_Kowalski}). In \cite[Proposition 4.1]{Console_Sergio_Fino_Evangelia}, they list the standard representatives of 7 isomorphism classes of 2-step nilpotent Lie algebras of dimension 6 (together with $\mathbb{R}^6$, abelian). Among their list, (5) and (8) correspond to $\textrm{heis}_3(\mathbb{R}[\epsilon]), \mathcal{N}_6$, respectively. To see this, one can just observe that $\textrm{heis}_3(\mathbb{R}[\epsilon])$ is of type $(2,4)$, while $\mathcal{N}_6$ is of type $(3,3)$. 

\begin{proof}[Proof of Theorem \ref{pac_up_to_dim_6_introduction}]
\textbf{Dimensions 3,4, and 5}. Recall, we may either work at the Lie group level or Lie algebra level by Definition \ref{def_of_partial_automatic_continuity_at_the_Lie_algebra_level} and Lemma \ref{reduction_to_the_Lie_algebras}. By Theorem \ref{temp_theorem1_intro}, which is independently shown in Section \ref{section_definitions}, $\textrm{heis}_3(\mathbb{R})$ and $\textrm{heis}_5(\mathbb{R})$ satisfy partial automatic continuity. By Theorem \ref{maintheorem_criterion}, $\mathcal{N}_5$ satisfies partial automatic continuity. 
By Proposition \ref{abelian_factors_2-step} (independently shown in Section \ref{section_abelian_factors}), all the 2-step nilpotent Lie algebras of dimensions 3,4, and 5 satisfy partial automatic continuity.

\textbf{Dimension 6.} By Proposition \ref{abelian_factors_2-step} and the results above, we only need to check $\textrm{heis}_3(\mathbb{R})\oplus \textrm{heis}_3(\mathbb{R})$, $\textrm{heis}_3(\mathbb{C})$, $\textrm{heis}_3(\mathbb{R}[\epsilon])$, and $\mathcal{N}_6$. The case $\textrm{heis}_3(\mathbb{C})$ is by the case (2) of Theorem \ref{maintheorem_over_C}. Theorem \ref{maintheorem_criterion} covers $\mathcal{N}_6$, as $[X_1, X_2]=Z_1$ and $[X_1,X_3]=-Z_2$ are linearly independent. The case $\textrm{heis}_3(\mathbb{R})\oplus \textrm{heis}_3(\mathbb{R})$ is covered by Corollary \ref{corollary_concatenation} (1) as it is a concatenation of admissible 2-step nilpotent Lie algebras. The case of $\textrm{heis}_3(\mathbb{R}[\epsilon])$ is covered by \cite[p.554, Corollary 5]{Levchuk}.
\end{proof}

\section{Proof of Theorem \ref{temp_theorem1_intro} and Case (3) of Theorem \ref{maintheorem_over_C}}\label{section_definitions}

The Iwasawa N-group is the N part of the Iwasawa decomposition $G=KAN$ of a simple Lie group $G$. The Iwasawa N-group $N$ of rank 1 simple Lie groups $SO_0(n,1), SU(n,1), Sp(n,1)$, and F II are, respectively, abelian, Heisenberg groups, quaternionic Heisenberg groups, and the 15-dimensional octonionic Heisenberg group (see \cite{Cowling_Dooley_Koranyi_Ricci}, \cite{Calin_Chang_Markina}). The first three arise as families. We first describe the corresponding Lie algebra.

The $(2n+1)$-dimensional Heisenberg algebra $\textrm{heis}_{2n+1}(\mathbb{R})$ is a $(2n+1)$-dimensional 2-step nilpotent Lie algebra with the following non-trivial Lie brackets: $[X_1, Y_1]=...=[X_n, Y_n]=Z$, where we declare that $X_1, Y_1, ..., X_n, Y_n, Z$ is a $\mathbb{R}$-basis of $\textrm{heis}_{2n+1}(\mathbb{R})$ (\cite[p.617, Example 1]{Eberlein}. The $(4n+3)$-dimensional quaternionic Heisenberg algebra $\textrm{q-heis}_{4n+3}(\mathbb{R})$ is a  $(4n+3)$-dimensional 2-step nilpotent Lie algebra with the following non-trivial Lie brackets: $[X_i, Y_i]=Z_1, [X_i, V_i]=Z_2, [X_i,W_i]=Z_3, [Y_i, V_i]=Z_3, [Y_i,W_i]=-Z_2, [V_i,W_i]=Z_1$, where we declare that $X_i, Y_i, V_i, W_i$, $i\in \{1,...,n\}$, and $Z_1, Z_2, Z_3$ are a $\mathbb{R}$-basis of $\textrm{q-heis}_{4n+3}(\mathbb{R})$ (\cite[p.617, Example 2]{Eberlein}). The 15-dimensional octonionic Heisenberg algebra $\textrm{o-heis}_{15}(\mathbb{R})$ has a $\mathbb{R}$-basis $X_1,...,X_8, Z_1, ..., Z_7$ with the Lie bracket relations given by \cite[p. 51, Table 3]{Calin_Chang_Markina}.

\begin{proof} [Proof of Theorem \ref{temp_theorem1_intro}.]

From the definitions, the 3-dimensional Heisenberg algebra, 7-dimensional quaternionic Heisenberg algebra, and 15-dimensional octonionic Heisenberg algebra are admissible 2-step nilpotent Lie algebras. For example, it follows from \cite[p. 51, Table 3]{Calin_Chang_Markina} that $[X_1, X_j]=-Z_{j-1}$ for $j\in \{2,..., 8\}$, so $\textrm{o-heis}_{15}(\mathbb{R})$ is admissible. The $(2n+1)$-dimensional Heisenberg algebras and $(4n+3)$-dimensional quaternionic Heisenberg algebras are concatenations of admissible 2-step nilpotent Lie algebras. Therefore, Corollary \ref{corollary_concatenation} (1) applies.
\end{proof}

\begin{proof}[Proof of Case (3) of Theorem \ref{maintheorem_over_C}]
    If $\mathcal{N}'$ denotes any of $\textrm{heis}_{2n+1}(\mathbb{R}), \textrm{q-heis}_{4n+3}(\mathbb{R}),$ and $\textrm{o-heis}_{15}(\mathbb{R})$, then the structure constants of the given basis are integers. Hence, the same basis works as a $\mathbb{C}$-basis of $\mathcal{N}=\mathcal{N}'\otimes_{\mathbb{R}}\mathbb{C}$. By Corollary \ref{corollary_concatenation} (2), the result follows.
\end{proof}

One might hope that Theorem \ref{temp_theorem1_intro} extends to any homomorphisms. We give a counterexample.

\begin{example}\label{discontinuous_homomorphism_preserving_v}
    Let $\mathcal{N}$ be the $(2n+1)$-dimensional Heisenberg algebra over $\mathbb{R}$, and let $\mathcal{N}=\mathcal{V}\oplus\mathcal{Z}$, where $\mathcal{V}=\mathbb{R}\textrm{-span}\{X_1,Y_1,...,X_n,Y_n\}$ and $\mathcal{Z}=\mathbb{R}\textrm{-span}\{Z\}$. Let $\lambda:\mathbb{R}\to \mathbb{R}$ be a discontinuous $\mathbb{Q}$-linear map. Define a $\mathbb{Q}$-linear map $f$ by setting $f(xX_1)=\lambda(x)X_1$ and $f=0$ on the $\mathbb{R}$-span of $Y_1, X_2,Y_2,...,X_n,Y_n,Z$.

    Then, $f$ is a discontinuous Lie ring homomorphism under which $\mathcal{V}$ is invariant, and $f|_{\mathcal{Z}}=0$. The map $f$ induces a discountinuous map $\mathcal{N}/\mathcal{Z}\to \mathcal{N}/\mathcal{Z}$. At the Lie group level, $F=\exp\circ f\circ \exp^{-1}$ is a discontinuous group homomorphism that induces a discontinuous map $N/Z\to N/Z$, where $N$ is the $(2n+1)$-dimensional Heisenberg group and $Z$ is its center. In particular, Theorem \ref{temp_theorem1_intro} does not extend to group homomorphisms. 
    
\end{example}

\begin{remark}
    The construction above works for any non-abelian nilpotent Lie algebra $\mathcal{N}$ over $\mathbb{R}$. Since $\mathcal{N}$ is nilpotent, there is a non-zero subspace $\mathcal{V}$ such that $\mathcal{N}=\mathcal{V}\oplus[\mathcal{N},\mathcal{N}]$. Let $X_1, ..., X_m\in \mathcal{V}$ be a basis of $\mathcal{V}$. We may assume that $X_1$ is not in the center. Let $\lambda:\mathbb{R}\to \mathbb{R}$ be a discontinuous $\mathbb{Q}$-linear map. Define a $\mathbb{Q}$-linear map $f$ by setting $f(xX_1)=\lambda(x)X_1$, $f(X_2)=...=f(X_m)=0$, and $f|_{[\mathcal{N},\mathcal{N}]}=0$
    
    Then, $f$ is a discontinuous Lie ring homomorphism such that $f(\mathcal{V})\subseteq \mathcal{V}$ and $f|_{[\mathcal{N}, \mathcal{N}]}=0$, and the same conclusion as Example \ref{discontinuous_homomorphism_preserving_v} holds.
\end{remark} 

\section{Proof of Theorem \ref{isomorphism_theorem_introduction}} \label{section_isomorphism_rigidity}

Let $G$ be a connected Lie group. Let $N$ be a simply connected 2-step nilpotent Lie group. First,

\begin{lemma}\label{rigidity_simply_connectedness} Let $G$ be a connected Lie group. Let $N$ be a simply connected 2-step nilpotent Lie group.
    If $F:G\to N$ is an abstract group isomorphism, then $G$ is a simply connected 2-step nilpotent Lie group.
\end{lemma}

\begin{proof}
    Let $\tilde{G}$ be a universal covering Lie group of $G$, and let $\pi:\tilde{G}\to G$ be the universal covering Lie group homomorphism. Let $\tilde{\mathfrak{g}}, \mathfrak{g}, \mathcal{N}$ be the Lie algebras of $\tilde{G}, G, N$, respectively. (Of course, $\tilde{\mathfrak{g}}$ and $\mathfrak{g}$ are isomorphic.)
    
    First, $G$ is a 2-step nilpotent group as an abstract group. For a connected Lie group $H$ with Lie algebra $\mathfrak{h}$, $[H,H]$ is a connected Lie subgroup whose Lie algebra is $[\mathfrak{h}, \mathfrak{h}]$ (\cite[pp.125--127]{Chevalley_Lie_Groups}). Since $[G,G]\subseteq Z(G)$, the center, it follows $[\mathfrak{g}, \mathfrak{g}]\subseteq \mathfrak{z}(\mathfrak{g})$, the center of $\mathfrak{g}$. Thus, $G$ is a 2-step nilpotent Lie group, and $\tilde{G}$ is a simply connected nilpotent Lie group.

    We will show $\pi$ is bijective. Let $f=\exp^{-1}\circ (F\circ \pi) \circ \exp: \tilde{\mathfrak{g}}\to \mathcal{N}$. Here $\exp$ denotes the exponential maps of both $\tilde{G}$ and $N$. By Theorem \ref{correspondence}, $f$ is a surjective Lie ring homomorphism. 
    \vskip\baselineskip
    \noindent
    \textbf{Claim.} $\ker f = \exp^{-1} (\ker \pi)$.
    \vskip\baselineskip
    \noindent
    If $x\in \exp^{-1} (\ker \pi)$, then $f(x)=\exp^{-1}\circ F\circ \pi(\exp(x))=0$. Thus, $x\in \ker f$. If $x\in \ker f$, then $F(\pi(\exp(x)))=e$. As $F$ is injective, $\pi(\exp(x))=e$, proving Claim.

    In particular, $\ker f$ is discrete as $\exp$ is a diffeomorphism. If $0\neq x\in \ker f$, then for any $n\in \mathbb{N}$, $nf(\frac{1}{n}x)=f(x)=0$ as $f$ is additive. Thus, $f(\frac{1}{n}x)=0$, so $\frac{1}{n}x\in \ker f$. For any $m\in \mathbb{Z}$, $f(\frac{m}{n}x)=mf(\frac{1}{n}x)=0$, so $\frac{m}{n}x\in \ker f$. This is a contradiction as $\ker f$ is discrete. Therefore, $\ker f=\{0\}$. By Claim, $\ker \pi =\{e\}$, so $\pi$ is bijective. Hence, $G$ is a simply connected 2-step nilpotent Lie group. 
\end{proof}

We also need the following fact.

\begin{lemma}\label{Hausdorff_formula}
For any $x,y\in \mathfrak{g}$, 
\begin{align*}
    \exp(x)\exp(y)=\exp\left(x+y+\frac{1}{2}[x,y]\right)=\exp(x+y)\exp\left(\frac{1}{2}[x,y]\right).
\end{align*}
Furthermore, we have $\exp([x,y])=[\exp(x), \exp(y)]$, where the bracket on the left-hand side is the Lie bracket, and the bracket on the right-hand side is the commutator $[X, Y]=XYX^{-1}Y^{-1}$ in $G$. 

\end{lemma}

This is the Baker--Campbell--Hausdorff formula. Since the commutator ideal is central, the exponential map preserves the brackets. 

\begin{lemma} \label{induced_diffeo_from_exp}
    Let $N$ be a simply connected 2-step nilpotent Lie group with Lie algebra $\mathcal{N}$. Let $Z, \mathcal{Z}$ be the center of $N, \mathcal{N}$, respectively. 

    Then, the exponential map $\exp:\mathcal{N}\to N$ is a diffeomorphism, $\exp(\mathcal{Z})=Z$, $\exp([\mathcal{N},\mathcal{N}])=[N,N]$, and it induces a diffeormophism $\overline{\exp}:\mathcal{N}/\mathcal{Z}\to G/Z$.
\end{lemma}
\begin{proof}
    For any simply connected nilpotent Lie group, $\exp$ is a diffeomorphism. By Lemma \ref{Hausdorff_formula}, we have $\exp([x,y])=[\exp(x),\exp(y)]$. This implies $\exp(\mathcal{Z})=Z$ and $\exp([\mathcal{N},\mathcal{N}])=[N,N]$. If $x-y\in \mathcal{Z}$, then $\exp(x)=\exp(y+x-y)=\exp(y)\exp(x-y)$ with $\exp(x-y)\in Z$. Thus, $\overline{\exp}(x+\mathcal{Z})=\exp(x)+Z$ is well defined. Similarly, $\overline{\exp^{-1}}(xZ)=\exp^{-1}(x)+\mathcal{Z}$ is a well-defined inverse of $\overline{\exp}$. As $Z$ is a closed subgroup, the quotient map $N\to N/Z$ is a smooth submersion, so $\overline{\exp}$ and $\overline{\exp^{-1}}$ are smooth. This completes the proof.
\end{proof}

We are ready to prove 
\begin{theorem} \label{theorem_isomorphism_abstract_and_Lie_group}
    Let $G$ be a connected Lie group, and let $N$ be a simply connected 2-step nilpotent Lie group with Lie algebra $\mathcal{N}$. Suppose $\mathcal{N}$ is a concatenation of admissible 2-step nilpotent Lie algebras as defined in Definition \ref{definition_concatenation}. Let $Z(G), Z(N)$ be the center of $G, N$.

    If $F:G\to N$ is an abstract group isomorphism, then $G$ is a simply connected 2-step nilpotent Lie group, and the induced map $\overline{F}:G/Z(G)\to N/Z(N)$ and $F|_{[G,G]}:[G,G]\to [N,N]$ are diffeomorphisms. 

    If in addition $[N,N]=Z(N)$, then there is a Lie group isomorphism from $G$ to $N$. 
\end{theorem}
\begin{proof}
    By Lemma \ref{rigidity_simply_connectedness}, $G$ is a simply connected 2-step nilpotent Lie group. Let $\tilde{\mathcal{N}}, \mathcal{N}$ be Lie algebras of $G, N$, respectively. Let $\tilde{\mathcal{Z}}, \mathcal{Z}$ be the center of $\tilde{\mathcal{N}}, \mathcal{N}$, respectively. Let $f=\exp^{-1}\circ F\circ \exp:\tilde{\mathcal{N}}\to \mathcal{N}$ be the corresponding Lie ring isomorphism from Theorem \ref{correspondence}. Then, Setting \ref{setting_from_arbitrary_to_concatenation} is satisfied with $\mathbb{F}=\mathbb{R}$. We use the notations in Notation \ref{notation__from_arbitrary_to_concatenation}. 

    We show $\tilde{\mathcal{N}}=\tilde{\mathcal{V}} \oplus \tilde{\mathcal{Z}}$. Let $x\in \tilde{\mathcal{N}}$. Then write $f(x)\in \mathcal{N}$ as $f(x) =\sum_{j=1}^Q a_j e_j +z$, where $a_j\in \mathbb{R}$ and $z\in \mathcal{Z}$. By the item 5 of Lemma \ref{lemma_from_arbitrary_to_concatenation}, $f(\sum_{j=1}^Q a_j \tilde{e}_j)\in \sum_{j=1}^Q a_j e_j +\mathcal{Z}$. We have $f(x-\sum_{j=1}^q a_j \tilde{e}_j)\in \mathcal{Z}$. Thus, $x-\sum_{j=1}^q a_j \tilde{e}_j \in f^{-1}(\mathcal{Z})=\tilde{\mathcal{Z}}$ by Lemma \ref{center_goes_to_center}. This shows $\tilde{\mathcal{N}}=\tilde{\mathcal{V}} + \tilde{\mathcal{Z}}$. If $x\in \tilde{\mathcal{V}}\cap \tilde{\mathcal{Z}}$, we can write $x=\sum_{j=1}^Q a_j \tilde{e}_j$. Then, $f(\sum_{j=1}^Q a_j \tilde{e}_j)= \sum_{j=1}^Q a_j e_j +z$, $z\in \mathcal{Z}$ by Lemma \ref{lemma_from_arbitrary_to_concatenation}. On the other hand, $f(x)\in \mathcal{Z}$, so $\sum_{j=1}^Q a_j e_j\in \mathcal{V}\cap \mathcal{Z}=0$. Thus, $a_j=0$ for any $j\in \{1,...,Q\}$, so $x=0$. This proves $\tilde{\mathcal{N}}=\tilde{\mathcal{V}} \oplus \tilde{\mathcal{Z}}$.
    
    Let $\overline{f}:\tilde{\mathcal{N}}/\tilde{\mathcal{Z}}\to \mathcal{N}/\mathcal{Z}$ be the induced map 
    $\overline{f}(x+\tilde{\mathcal{Z}})=f(x)+\mathcal{Z}$. Since $f$ is bijective, so is $\overline{f}$. We show $\overline{f}$ is a $\mathbb{R}$-linear isomorphism. Any element of $\tilde{\mathcal{N}}/\tilde{\mathcal{Z}}$ can be written as $x+\tilde{\mathcal{Z}}, x\in \tilde{\mathcal{V}}$ by the argument above. Let $t\in \mathbb{R}$. Then, by using the item 5 of Lemma \ref{lemma_from_arbitrary_to_concatenation}, $\overline{f}(tx+\tilde{\mathcal{Z}})=f(tx)+\mathcal{Z}=tf(x)+\mathcal{Z}=t\overline{f}(x+\tilde{\mathcal{Z}})$. Since $\overline{f}$ is bijective, this shows $\overline{f}$ is a $\mathbb{R}$-linear isomorphism. Note, $\overline{F}(\overline{\exp}(x+\tilde{\mathcal{Z}}))=F(\exp(x))Z(N)=\exp(f(x))Z(N)=\overline{\exp}(\overline{f}(x+\tilde{\mathcal{Z}})$. Thus, by Lemma \ref{induced_diffeo_from_exp}, $\overline{F}$ is a 
    diffeomorphism.  

    Next, we show $f|_{[\tilde{\mathcal{N}}, \tilde{\mathcal{N}}]}: [\tilde{\mathcal{N}}, \tilde{\mathcal{N}}] \to [\mathcal{N},\mathcal{N}]$ is $\mathbb{R}$-linear. Since $\tilde{\mathcal{N}}=\tilde{\mathcal{V}} \oplus \tilde{\mathcal{Z}}$, we have $[\tilde{\mathcal{N}}, \tilde{\mathcal{N}}]=[\tilde{\mathcal{V}},\tilde{\mathcal{V}}]$. Any $z\in [\tilde{\mathcal{N}}, \tilde{\mathcal{N}}]$ has the form $z=\sum_{j=1}^l [x_j,y_j]$, $x_j,y_j\in \tilde{\mathcal{V}}$. For any $t\in \mathbb{R}$, we have $f(tz)=\sum_{j=1}^l [f(tx_j), f(y_j)]=t\sum_{j=1}^l[f(x_j),f(y_j)]=tf(z)$. Thus, $f|_{[\tilde{\mathcal{N}}, \tilde{\mathcal{N}}]}$ is $\mathbb{R}$-linear. Since it has the inverse $f^{-1}|_{[\mathcal{N},\mathcal{N}]}$, $f|_{[\tilde{\mathcal{N}}, \tilde{\mathcal{N}}]}$ is a $\mathbb{R}$-linear isomorphism. By Lemma \ref{induced_diffeo_from_exp}, $F|_{[G,G]}$ is a diffeormopshim.

    Finally, we show that there is a Lie group isomorphism from $G$ to $N$ assuming $[N,N]=Z(N)$. As $[N,N]=Z(N)$ and $F$ is a group isomorphism, we have $[G,G]=Z(G)$. By the results above, $\dim G=\dim G/Z(G) + \dim [G,G]=\dim N/Z(N) + \dim [N,N] =\dim N$. Let $\piv:\mathcal{N}\to \mathcal{V}$ and $\piz:\mathcal{N}\to \mathcal{Z}$ be the projections such that $\piv(X+Z)=X, \piv(X+Z)=Z$, where $X\in \mathcal{V}$ and $Z\in \mathcal{Z}$.
    
    Let us define $g:\tilde{\mathcal{N}}\to \mathcal{N}$ by $g(X+Z)=\piv\circ f(X)+f(Z)$ for $X\in \tilde{\mathcal{V}}, Z\in \tilde{\mathcal{Z}}$. We show $g$ is a Lie algebra isomorphism. First, observe that $g$ is additive. For any $x,y\in \tilde{\mathcal{N}}$, $g([x,y])=f([x,y])=[f(x),f(y)]=[g(x),g(y)]$. Thus, $g$ is a Lie ring homomorphism. 

    Next, we show $g$ is injective. As $f=g$ on $\tilde{\mathcal{Z}}$, $g|_{\tilde{\mathcal{Z}}}$ is injective. Let $X=\sum_{j=1}^Q a_j \tilde{e}_j\in \tilde{\mathcal{V}}$ and $g(X)=0$. By Lemma \ref{lemma_from_arbitrary_to_concatenation}, $f(\sum_{j=1}^Q a_j \tilde{e}_j)=\sum_{j=1}^Q a_j e_j +z$ for some $z\in \mathcal{Z}$. Thus, $0=\piv(f(X))=\piv (\sum_{j=1}^Q a_j e_j+z)=\sum_{j=1}^Q a_j e_j$. As $e_1,...,e_Q$ is linearly independent, $a_j=0$ for any $j\in \{1,...,Q\}$. Thus, $X=0$ and $g$ is injective. We show $g$ is $\mathbb{R}$-linear. For any $X\in \tilde{\mathcal{V}}$ write $X=\sum_{j=1}^Q a_j \tilde{e}_j\in \tilde{\mathcal{V}}$. For any $t\in \mathbb{R}$, $g(tX)=\sum_{j=1}^Q ta_j e_j =tg(X)$ by Lemma \ref{lemma_from_arbitrary_to_concatenation} and $g(\tilde{\mathcal{V}})\subseteq \mathcal{V}$. As $f(z)=g(z)$ for $z\in [\tilde{\mathcal{N}}, \tilde{\mathcal{N}}]=\tilde{\mathcal{Z}}$, $f|_{\tilde{\mathcal{Z}}}=g|_{\tilde{\mathcal{Z}}}$ is $\mathbb{R}$-linear. Thus, $g$ is $\mathbb{R}$-linear. 

    Since $\dim G=\dim N$, we have $\dim \mathcal{N}=\dim \tilde{\mathcal{N}}$. Since $g$ is an injective $\mathbb{R}$-linear map between finite-dimensional vector spaces of the same dimension, $g$ is bijective. Therefore, $g$ is a Lie algebra isomorphism. As $G$ and $N$ are simply connected, there is a Lie group isomorphism between $G$ and $N$. This completes the proof.
\end{proof}

\section{Semidirect Product Structure of Abstract Group Automorphisms}\label{section_semidirect_product}

Let $\mathbb{F}$ be any field. This section is independent of the other sections. Any 2-step nilpotent Lie algebra $\mathcal{N}$ admits a Carnot grading: if $\mathcal{Z}=[\mathcal{N}, \mathcal{N}]$ and $\mathcal{V}$ is any $\mathbb{F}$-subspace of $\mathcal{N}$ with $\mathcal{N}=\mathcal{V}\oplus \mathcal{Z}$, then $[\mathcal{V}, \mathcal{V}]=\mathcal{Z}$. This is also a Carnot grading as a Lie ring, or in other words, as a $\mathbb{Z}$-algebra (see \cite[Section 3.1]{Cornulier} for the definition), where a $\mathbb{Z}$-bilinear law is given by the Lie bracket. By \cite[Corollary 3.6]{Cornulier}, any Lie ring automorphism $f$ is a product $f=\mu \circ \overline{f}$, where $\mu$ is a Lie ring automorphism inducing the identity on $\mathcal{N}/\mathcal{Z}$ and $\overline{f}$ is a Lie ring automorphism such that $f(\mathcal{V})=\mathcal{V}$. In particular, if one only cares about 2-step nilpotent Lie algebras without an abelian factor, then one needs Proposition \ref{semidirect_product_for_2-step_nilpotent_Lie_algebras_over_fields} below in the case of $\mathcal{Z}=[\mathcal{N}, \mathcal{N}]$, which is a special case of \cite[Corollary 3.6]{Cornulier}.

However, we allow an abelian factor in Theorem \ref{maintheorem_criterion}, Theorem \ref{pac_up_to_dim_6_introduction}, and the case (2) of Theorem \ref{maintheorem_over_C}. For these results, we prove the following proposition.

\begin{proposition} \label{semidirect_product_for_2-step_nilpotent_Lie_algebras_over_fields} Let $\mathbb{F}$ be any field.
    Let $\mathcal{N}$ be a finite-dimensional 2-step nilpotent Lie algebra over $\mathbb{F}$. Let $\mathcal{Z}$ be $[\mathcal{N}, \mathcal{N}]$ or $\mathfrak{z}(\mathcal{N})$, the center of $\mathcal{N}$. Let $\mathcal{V}$ be a subspace of $\mathcal{N}$ such that $\mathcal{N}=\mathcal{V}\oplus \mathcal{Z}$. Let $G$ be the group of Lie ring isomorphisms of $\mathcal{N}$ as a Lie ring. Let $H=\{f\in G\;|\; \textrm{for any $x\in \mathcal{N}$, }f(x)-x\in \mathcal{Z}\}$ and $K=\{f\in G\;|\;f(\mathcal{V})=\mathcal{V}\}$.
    Then, $H$ is a normal subgroup of $G$, $K$ is a subgroup of $G$, and $G=HK$. In the case $\mathcal{Z}=[\mathcal{N}, \mathcal{N}]$, we have $H\cap K=\{id\}$, and thus $G \cong H\rtimes K$.
\end{proposition}
\begin{remark}
    The group $H$ is the group of central automorphisms if $\mathcal{Z}=\mathfrak{z}(\mathcal{N})$.
\end{remark}

\begin{proof}
    Define $\pi_{\mathcal{V}}: \mathcal{N}\to \mathcal{V}$ and $\pi_{\mathcal{Z}}: \mathcal{N}\to \mathcal{Z}$ by $\piv(X+Z)=X$ and $\piz(X+Z)=Z$, where $X\in \mathcal{V}$ and $Z\in \mathcal{Z}$. Let $\iota_{\mathcal{V}}:\mathcal{V}\to \mathcal{N}$ and $\iota_{\mathcal{Z}}:\mathcal{Z}\to \mathcal{N}$, be the inclusions. The key is that a Lie ring automorphism $f$ satisfies $f(\mathcal{Z})=\mathcal{Z}$. 
    
    By definition, $K$ is a subgroup. Let us show that $H$ is a subgroup. For $f\in H$, $f^{-1}(\mathcal{Z})=\mathcal{Z}$ because $f$ preserves the Lie bracket. Hence, $x-f^{-1}(x)=f^{-1}(f(x)-x)\in f^{-1}(\mathcal{Z})=\mathcal{Z}$. Thus, $f^{-1}\in H$. If $f_1, f_2\in H$, then $f_1\circ f_2(x)-x=f_1(f_2(x))-f_2(x)+f_2(x)-x\in \mathcal{Z}$. Thus, $f_1\circ f_2\in H$. This shows that $H$ is a subgroup. 
    
    Let us show that $H$ is a normal subgroup. Let $f\in H$ and $g\in G$. Then, $g\circ f\circ g^{-1}(x)-x=g(f(g^{-1}(x))-g^{-1}(x))\in \mathcal{Z}$ since $g(\mathcal{Z})=\mathcal{Z}$. Thus, $H$ is a normal subgroup. 
    
    Next, let us show $G=HK$. By taking the inverse, it suffices to show that $G=KH$. Let $f\in G$. Let us write $h=f|_{\mathcal{Z}}:\mathcal{Z}\to \mathcal{Z}$. Since $f$ is a Lie ring automorphism, $h$ is bijective. Let $\mu(X+Z)=h^{-1}(\piz\circ f(X))$, where $X\in \mathcal{V}$ and $Z\in \mathcal{Z}$. Let $f_1=\textrm{id}_{\mathcal{N}}+\mu$. Since $\mu|_{[\mathcal{N},\mathcal{N}]}=0$ and $\mu^2=0$, it follows from Lemma \ref{characterization_central_automorphism} that $f_1$ is a central automorphism. Thus, $f_1\in H$. Let $f_2(X+Z)=\pi_{\mathcal{V}}(f(X))+f(Z)$, where $X\in \mathcal{V}$ and $Z\in \mathcal{Z}$,  so that $f_2(\mathcal{V})\subseteq \mathcal{V}$.
    
    We show $f=f_2\circ f_1$. For any $Z\in \mathcal{Z}$, we have $f_2(f_1(Z))=f(Z)$. For any $X\in \mathcal{V}$, we have $f_2(f_1(X))=f_2(X+h^{-1}(\piz\circ f(X))=\piv(f(X))+\piz(f(X))=f(X)$. Thus, $f=f_2\circ f_1$. In particular, $f_2=f\circ f_1^{-1}$ is a Lie ring automorphism. To show  $f_2(\mathcal{V})=\mathcal{V}$, let $Y\in \mathcal{V}$. Set $X=f^{-1}(Y)-\piz \circ f^{-1}(Y)\in \mathcal{V}$. Then, noting $f=f_2$ modulo $\mathcal{Z}$, we have $f_2(X)\in f_2(f^{-1}(Y))+\mathcal{Z}=f(f^{-1}(Y))+\mathcal{Z}=Y+\mathcal{Z}$. Thus, $f_2(X)-Y\in \mathcal{V}\cap \mathcal{Z}=\{0\}$, so $f_2(X)=Y$. Hence, $f_2\in K$.
    
    Finally, suppose $\mathcal{Z}=[\mathcal{N}, \mathcal{N}]$. We show $H\cap K=\{\textrm{id}\}$. If $f\in H\cap K$, for any $X\in \mathcal{V}$, $f(X)-X\in \mathcal{V}\cap \mathcal{Z}=\{0\}$. Thus, $f|_{\mathcal{V}}$ is the identity. Since $\mathcal{Z}=[\mathcal{V}, \mathcal{V}]$, $f|_{\mathcal{Z}}$ is also the identity. This completes the proof.
\end{proof}

\begin{remark}\label{remark_on_semidirect_product_structure_at_group_level}
    This proposition generalizes the smooth case \cite[Proposition 3.4.2]{Eberlein_semidirect_product} in the case of $\mathbb{F}=\mathbb{R}$. By using Theorem \ref{correspondence}, we see that the group of abstract group automorphisms of a simply connected 2-step nilpotent Lie group has a semidirect product structure as well.
\end{remark}

\begin{remark} Let $\mathbb{F}=\mathbb{R}$. 
    Let $\mathcal{N}=\mathcal{V}\oplus [\mathcal{N},\mathcal{N}]$ and $p=\dim [\mathcal{N},\mathcal{N}]$ and $q=\dim \mathcal{V}$. The group $H$ is isomorphic to $\{\lambda:\mathbb{R}^q\to \mathbb{R}^p\;|\;\textrm{$\lambda$ is $\mathbb{Q}$-linear}\}$. For a generic algebra, the center of $\mathcal{N}$ coincides with $[\mathcal{N}, \mathcal{N}]$. When the condition in Theorem \ref{maintheorem_criterion} is satisfied, the group $K$ is isomorphic to $\{f\in \textrm{Aut}(\mathcal{N})\;|\;f(\mathcal{V})=\mathcal{V}\}$ by Theorem \ref{theorem_concatenation}, where $\textrm{Aut}(\mathcal{N})$ denotes the group of Lie algebra automorphisms. We always have $\dim K\geq 1$ because $\mathcal{N}$ always admits a symmetric derivation. The group $K$ corresponds to $\overline{H}$ in \cite[Proposition 3.4.2]{Eberlein_semidirect_product}, see Remark below Proposition 3.4.2. Suppose $(p,q)$ is not in Table in \cite[p.87]{Eberlein_semidirect_product}. By \cite[Proposition 3.4.2]{Eberlein_semidirect_product} and ``Reduction to $SL(q,\mathbb{R})$'' in \cite[p.82]{Eberlein_semidirect_product}, $K$ is generically one-dimensional, corresponding to the $(1,2)$-derivation, i.e., a derivation that is $\textrm{id}_{\mathcal{V}}$ on $\mathcal{V}$ and $2\textrm{id}_{[\mathcal{N},\mathcal{N}]}$ on $[\mathcal{N},\mathcal{N}]$. Hence, the group of abstract group automorphisms of a generic simply connected 2-step nilpotent Lie group with $\dim [N,N]\geq \frac{\dim N-1}{2}$ is the smallest.
\end{remark}

\section{Abelian Factors} \label{section_abelian_factors}

Let $\mathbb{F}$ be a field. The result in this section only depends on Section \ref{section_semidirect_product}. For any 2-step nilpotent Lie algebra $\mathcal{N}$ over $\mathbb{F}$, there are ideals $\mathcal{N}_0$ and $\mathcal{A}$ such that $\mathcal{N}=\mathcal{N}_0\oplus \mathcal{A}$, $\mathcal{A}$ is abelian, and $[\mathcal{N}_0, \mathcal{N}_0]=\mathfrak{z}(\mathcal{N}_0)$, the center of $\mathcal{N}_0$. We call $\mathcal{A}$ an \textit{abelian factor} of $\mathcal{N}$ (see \cite[Section 1.3]{Eberlein_moduli_2-step}).

\begin{proposition}\label{abelian_factors_2-step} Let $\mathbb{F}$ be a field. Let $\mathcal{N}$ be a 2-step nilpotent Lie algebra over $\mathbb{F}$. Suppose $\mathcal{N}=\mathcal{N}_0\oplus \mathcal{A}$ (a direct sum of ideals) with $\mathcal{A}$ the abelian factor and $[\mathcal{N}_0, \mathcal{N}_0]=\mathfrak{z}(\mathcal{N}_0)$, the center of $\mathcal{N}_0$. Then, any Lie ring automorphism $f$ has the form $f=\mu\circ g$, where $\mu$ is a central automorphism and $g$ is a Lie ring automorphism such that $g(\mathcal{N}_0)=\mathcal{N}_0$ and $g(Y)=Y$ for $Y\in \mathcal{A}$. In particular, 
\begin{enumerate}
    \item If $\mathbb{F}=\mathbb{R}$ and $\mathcal{N}_0$ satisfies partial automatic continuity, then so does $\mathcal{N}$.
    \item If any Lie ring automorphism of $\mathcal{N}_0$ is a product of a central automorphism and a $\varphi$-linear Lie ring automorphism for some $\varphi\in \textrm{Aut}(\mathbb{F})$, then the same holds for any Lie ring automorphism of $\mathcal{N}$. 
\end{enumerate}

\end{proposition}
\begin{proof} 
    Let $\mathcal{V}$ be a subspace so that $\mathcal{N}_0=\mathcal{V}\oplus \mathfrak{z}(\mathcal{N}_0)$. Let $\mathfrak{z}(\mathcal{N})=\mathfrak{z}(\mathcal{N}_0)\oplus \mathcal{A}$ be the center of $\mathcal{N}$. Let $\mathcal{Z}=\mathfrak{z}(\mathcal{N})$ so that  $\mathcal{N}=\mathcal{V}\oplus \mathcal{Z}$. By Proposition \ref{semidirect_product_for_2-step_nilpotent_Lie_algebras_over_fields}, there exists a central automorphism $\lambda$ and a Lie ring automorphism $g$ with $g(\mathcal{V})=\mathcal{V}$ such that $\lambda\circ g=f$. Since $g$ is a Lie ring automorphism, it maps the commutator ideal $[\mathcal{N}, \mathcal{N}]=[\mathcal{N}_0, \mathcal{N}_0]$ onto the commutator ideal. Since $\mathcal{N}_0=\mathcal{V}\oplus [\mathcal{N}_0, \mathcal{N}_0]$, we have $g(\mathcal{N}_0)=\mathcal{N}_0$.
    
    Next, set $g_1(X+Y)=g(X)+Y$ and $g_2(X+Y)=X+g(Y)$ for $X\in \mathcal{N}_0$ and $Y\in \mathcal{A}$. Then, $g_2\circ g_1=g$. Note, $g_1$ and $g_2=g\circ g_1^{-1}$ are Lie ring automorphisms since $g_1([X_1+Y_1, X_2+Y_2])=g([X_1,X_2])=[g(X_1),g(X_2)]=[g_1(X_1+Y_1), g_1(X_2+Y_2)]$, where $X_1,X_2 \in \mathcal{N}_0$ and $Y_1, Y_2\in \mathcal{A}$.
    
    Thus, $f=\lambda\circ g_2\circ g_1$, $g_1(\mathcal{N}_0)=\mathcal{N}_0$, and $g_1(Y)=Y$ for any $Y\in \mathcal{A}$, and $g_2$ is a central automorphism.

    The item 1 immediately follows. To see the item 2, say $f$ is a Lie ring automorphism of $\mathcal{N}$ and $f=\mu \circ g$ as above. Then $g|_{\mathcal{N}_0}=\lambda\circ \overline{g}$, where $\lambda$ is a central automorphism and $\overline{g}$ is a $\varphi$-linear Lie ring automorphism of $\mathcal{N}_0$. Let $e_1,...,e_m$ be a basis of $\mathcal{A}$ and let $\overline{\varphi}(\sum_{j=1}^m a_je_j)=\sum_{j=1}^m\varphi(a_j)e_j$. 
    
    Now extend $\overline{g}$ to $\tilde{g}:\mathcal{N}\to \mathcal{N}$ by setting $\tilde{g}(x)=\overline{g}(x)$ for any $x\in \mathcal{N}_0$ and $\tilde{g}(x)=\overline{\varphi}(x)$ for any $x\in \mathcal{A}$, and define $h:\mathcal{N}\to \mathcal{N}$ by $h=\textrm{Id}$ on $\mathcal{N}_0$ and $h=\overline{\varphi}^{-1}$ on $\mathcal{A}$. Extend $\lambda$ to  $\tilde{\lambda}:\mathcal{N}\to \mathcal{N}$ by setting $\tilde{\lambda}=\textrm{Id}$ on $\mathcal{A}$. Then, $f=\mu \circ \tilde{\lambda} \circ h \circ \tilde{g}$, and $\tilde{\lambda}, h$ are central automorphisms. Since $\tilde{g}$ is $\varphi$-linear Lie ring automorphism of $\mathcal{N}$, the result follows.
\end{proof}

\bibliographystyle{alpha}
\bibliography{main}{}

@article {Car1930,
    AUTHOR = {Cartan, Elie},
     TITLE = {Sur les repr\'{e}sentations lin\'{e}aires des groupes clos},
   JOURNAL = {Comment. Math. Helv.},
  FJOURNAL = {Commentarii Mathematici Helvetici},
    VOLUME = {2},
      YEAR = {1930},
    NUMBER = {1},
     PAGES = {269--283},
      ISSN = {0010-2571},
   MRCLASS = {DML},
  MRNUMBER = {1509418},
       DOI = {10.1007/BF01214464},
       URL = {https://doi.org/10.1007/BF01214464},
}

@article {vanderwaerden1933,
    AUTHOR = {van der Waerden, B. L.},
     TITLE = {Stetigkeitss\"{a}tze f\"{u}r halbeinfache {L}iesche {G}ruppen},
   JOURNAL = {Math. Z.},
  FJOURNAL = {Mathematische Zeitschrift},
    VOLUME = {36},
      YEAR = {1933},
    NUMBER = {1},
     PAGES = {780--786},
      ISSN = {0025-5874},
   MRCLASS = {DML},
  MRNUMBER = {1545369},
       DOI = {10.1007/BF01188647},
       URL = {https://doi.org/10.1007/BF01188647},
}

@book {Segal,
    AUTHOR = {Segal, Daniel},
     TITLE = {Polycyclic groups},
    SERIES = {Cambridge Tracts in Mathematics},
    VOLUME = {82},
 PUBLISHER = {Cambridge University Press, Cambridge},
      YEAR = {1983},
     PAGES = {xiv+289},
      ISBN = {0-521-24146-4},
   MRCLASS = {20-02 (20F16)},
  MRNUMBER = {713786},
MRREVIEWER = {John F. Bowers},
       DOI = {10.1017/CBO9780511565953},
       URL = {https://doi.org/10.1017/CBO9780511565953},
}

@article {Eberlein,
    AUTHOR = {Eberlein, Patrick},
     TITLE = {Geometry of {$2$}-step nilpotent groups with a left invariant
              metric},
   JOURNAL = {Ann. Sci. \'{E}cole Norm. Sup. (4)},
  FJOURNAL = {Annales Scientifiques de l'\'{E}cole Normale Sup\'{e}rieure.
              Quatri\`eme S\'{e}rie},
    VOLUME = {27},
      YEAR = {1994},
    NUMBER = {5},
     PAGES = {611--660},
      ISSN = {0012-9593},
   MRCLASS = {53C22 (53C20 58F17)},
  MRNUMBER = {1296558},
MRREVIEWER = {Carolyn\ Gordon},
       URL = {http://www.numdam.org/item?id=ASENS_1994_4_27_5_611_0},
}

@incollection {Eberlein_semidirect_product,
    AUTHOR = {Eberlein, Patrick},
     TITLE = {Geometry of 2-step nilpotent {L}ie groups},
 BOOKTITLE = {Modern dynamical systems and applications},
     PAGES = {67--101},
 PUBLISHER = {Cambridge Univ. Press, Cambridge},
      YEAR = {2004},
      ISBN = {0-521-84073-2},
   MRCLASS = {53C30 (22E25 37D40)},
  MRNUMBER = {2090766},
MRREVIEWER = {Jorge\ Lauret},
}

@article {Calin_Chang_Markina,
    AUTHOR = {Calin, Ovidiu and Chang, Der-Chen and Markina, Irina},
     TITLE = {Geometric analysis on {$H$}-type groups related to division
              algebras},
   JOURNAL = {Math. Nachr.},
  FJOURNAL = {Mathematische Nachrichten},
    VOLUME = {282},
      YEAR = {2009},
    NUMBER = {1},
     PAGES = {44--68},
      ISSN = {0025-584X,1522-2616},
   MRCLASS = {53C17 (35H20)},
  MRNUMBER = {2473130},
MRREVIEWER = {Thomas\ Bieske},
       DOI = {10.1002/mana.200710721},
       URL = {https://doi.org/10.1002/mana.200710721},
}

@article {Braun_Hofmann_Kramer,
    AUTHOR = {Braun, O. and Hofmann, Karl H. and Kramer, L.},
     TITLE = {Automatic continuity of abstract homomorphisms between locally
              compact and {P}olish groups},
   JOURNAL = {Transform. Groups},
  FJOURNAL = {Transformation Groups},
    VOLUME = {25},
      YEAR = {2020},
    NUMBER = {1},
     PAGES = {1--32},
      ISSN = {1083-4362,1531-586X},
   MRCLASS = {22A05 (20K45 22A10 22C05 22D05 22E46)},
  MRNUMBER = {4070101},
MRREVIEWER = {Elena\ Mart\'{\i}n Peinador},
       DOI = {10.1007/s00031-019-09537-4},
       URL = {https://doi.org/10.1007/s00031-019-09537-4},
}

@article {Kallman_McLinden,
    AUTHOR = {Kallman, Robert R. and McLinden, Alexander P.},
     TITLE = {The {P}oincar\'{e} and related groups are algebraically
              determined {P}olish groups},
   JOURNAL = {Collect. Math.},
  FJOURNAL = {Universitat de Barcelona. Collectanea Mathematica},
    VOLUME = {61},
      YEAR = {2010},
    NUMBER = {3},
     PAGES = {337--352},
      ISSN = {0010-0757,2038-4815},
   MRCLASS = {22E43 (22A99 22E15)},
  MRNUMBER = {2732376},
MRREVIEWER = {Julien\ Melleray},
       DOI = {10.1007/BF03191237},
       URL = {https://doi.org/10.1007/BF03191237},
}

@article {Myasnikov_Sohrabi,
    AUTHOR = {Myasnikov, Alexei G. and Sohrabi, Mahmood},
     TITLE = {Groups elementarily equivalent to a free nilpotent group of
              finite rank},
   JOURNAL = {Ann. Pure Appl. Logic},
  FJOURNAL = {Annals of Pure and Applied Logic},
    VOLUME = {162},
      YEAR = {2011},
    NUMBER = {11},
     PAGES = {916--933},
      ISSN = {0168-0072,1873-2461},
   MRCLASS = {20A15 (03C60 20F18)},
  MRNUMBER = {2817564},
MRREVIEWER = {O.\ V.\ Belegradek},
       DOI = {10.1016/j.apal.2011.04.003},
       URL = {https://doi.org/10.1016/j.apal.2011.04.003},
}

@article {Gibbs,
    AUTHOR = {Gibbs, John A.},
     TITLE = {Automorphisms of certain unipotent groups},
   JOURNAL = {J. Algebra},
  FJOURNAL = {Journal of Algebra},
    VOLUME = {14},
      YEAR = {1970},
     PAGES = {203--228},
      ISSN = {0021-8693},
   MRCLASS = {20.22},
  MRNUMBER = {266985},
MRREVIEWER = {J.\ G.\ Horne},
       DOI = {10.1016/0021-8693(70)90123-7},
       URL = {https://doi.org/10.1016/0021-8693(70)90123-7},
}

@article {Borel_and_Tits,
    AUTHOR = {Borel, Armand and Tits, Jacques},
     TITLE = {Homomorphismes ``abstraits'' de groupes alg\'{e}briques
              simples},
   JOURNAL = {Ann. of Math. (2)},
  FJOURNAL = {Annals of Mathematics. Second Series},
    VOLUME = {97},
      YEAR = {1973},
     PAGES = {499--571},
      ISSN = {0003-486X},
   MRCLASS = {20G15 (10C30 14L15 20F55)},
  MRNUMBER = {316587},
MRREVIEWER = {J.\ Dieudonn\'{e}},
       DOI = {10.2307/1970833},
       URL = {https://doi.org/10.2307/1970833},
}

@article {Cowling_Dooley_Koranyi_Ricci,
    AUTHOR = {Cowling, Michael and Dooley, Anthony H. and Kor\'{a}nyi, Adam
              and Ricci, Fulvio},
     TITLE = {{$H$}-type groups and {I}wasawa decompositions},
   JOURNAL = {Adv. Math.},
  FJOURNAL = {Advances in Mathematics},
    VOLUME = {87},
      YEAR = {1991},
    NUMBER = {1},
     PAGES = {1--41},
      ISSN = {0001-8708,1090-2082},
   MRCLASS = {22E30 (22E46 22E60)},
  MRNUMBER = {1102963},
MRREVIEWER = {Detlef\ H.\ M\"{u}ller},
       DOI = {10.1016/0001-8708(91)90060-K},
       URL = {https://doi.org/10.1016/0001-8708(91)90060-K},
}

@article {Freudenthal_1941,
    AUTHOR = {Freudenthal, Hans},
     TITLE = {Die {T}opologie der {L}ieschen {G}ruppen als algebraisches
              {P}h\"{a}nomen. {I}},
   JOURNAL = {Ann. of Math. (2)},
  FJOURNAL = {Annals of Mathematics. Second Series},
    VOLUME = {42},
      YEAR = {1941},
     PAGES = {1051--1074},
      ISSN = {0003-486X},
   MRCLASS = {20.0X},
  MRNUMBER = {5740},
MRREVIEWER = {O.\ F. G. Schilling},
       DOI = {10.2307/1970461},
       URL = {https://doi.org/10.2307/1970461},
}

@article {Khor,
    AUTHOR = {Khor, Hoe Peng},
     TITLE = {The automorphisms of the unipotent radical of certain
              parabolic subgroups of {${\rm GL}(1+l,K)$}},
   JOURNAL = {J. Algebra},
  FJOURNAL = {Journal of Algebra},
    VOLUME = {96},
      YEAR = {1985},
    NUMBER = {1},
     PAGES = {54--77},
      ISSN = {0021-8693},
   MRCLASS = {20E36 (20G40)},
  MRNUMBER = {808841},
MRREVIEWER = {David\ B.\ Surowski},
       DOI = {10.1016/0021-8693(85)90039-0},
       URL = {https://doi.org/10.1016/0021-8693(85)90039-0},
}

@article {Im_symplectic,
    AUTHOR = {Im, Bokhee},
     TITLE = {The automorphisms of the unipotent radical of certain
              parabolic subgroups of {${\rm Sp}_{2l}(K)$}},
   JOURNAL = {Honam Math. J.},
  FJOURNAL = {Honam Mathematical Journal},
    VOLUME = {9},
      YEAR = {1987},
    NUMBER = {1},
     PAGES = {7--17},
      ISSN = {1225-293X},
   MRCLASS = {20G15},
  MRNUMBER = {919068},
MRREVIEWER = {R.\ W.\ Carter},
}

@article {Im_typeD,
    AUTHOR = {Im, Bokhee},
     TITLE = {On automorphisms of certain unipotent subgroups of {C}hevalley
              groups of type {$D^*_l$}},
   JOURNAL = {Bull. Korean Math. Soc.},
  FJOURNAL = {Bulletin of the Korean Mathematical Society},
    VOLUME = {28},
      YEAR = {1991},
    NUMBER = {1},
     PAGES = {23--32},
      ISSN = {1015-8634,2234-3016},
   MRCLASS = {20G15 (20E36)},
  MRNUMBER = {1107148},
MRREVIEWER = {Jean\ Michel},
}

@book {Im_dissertation,
    AUTHOR = {Im, Bokhee},
     TITLE = {The automorphisms of unipotent radicals of certain parabolic
              subgroups of symplectic groups},
      NOTE = {Thesis (Ph.D.)--University of Notre Dame},
 PUBLISHER = {ProQuest LLC, Ann Arbor, MI},
      YEAR = {1988},
     PAGES = {72},
   MRCLASS = {99-05},
  MRNUMBER = {2636313},
       URL =
              {http://gateway.proquest.com/openurl?url_ver=Z39.88-2004&rft_val_fmt=info:ofi/fmt:kev:mtx:dissertation&res_dat=xri:pqdiss&rft_dat=xri:pqdiss:8811934},
}

@book {Amayo_Stewart,
    AUTHOR = {Amayo, Ralph K. and Stewart, Ian},
     TITLE = {Infinite-dimensional {L}ie algebras},
 PUBLISHER = {Noordhoff International Publishing, Leiden},
      YEAR = {1974},
     PAGES = {xi+425},
   MRCLASS = {17B65},
  MRNUMBER = {396708},
MRREVIEWER = {Robert\ E.\ Beck},
}

@article {Malcev_Nilpotent_Torsion_Free_Groups,
    AUTHOR = {Malcev, A. I.},
     TITLE = {Nilpotent torsion-free groups},
   JOURNAL = {Izv. Akad. Nauk SSSR Ser. Mat.},
  FJOURNAL = {Izvestiya Akademii Nauk SSSR. Seriya Matematicheskaya},
    VOLUME = {13},
      YEAR = {1949},
     PAGES = {201--212},
      ISSN = {0373-2436},
   MRCLASS = {20.0X},
  MRNUMBER = {28843},
MRREVIEWER = {I.\ Kaplansky},
}

@incollection {Tits_Lie,
    AUTHOR = {Tits, J.},
     TITLE = {Homorphismes ``abstraits'' de groupes de {L}ie},
 BOOKTITLE = {Symposia {M}athematica, {V}ol. {XIII} ({C}onvegno di {G}ruppi
              {A}beliani \& {C}onvegno di {G}ruppi e loro
              {R}appresentazioni, {INDAM}, {R}ome, 1972)},
     PAGES = {479--499},
 PUBLISHER = {Academic Press, London-New York},
      YEAR = {1974},
   MRCLASS = {22E15},
  MRNUMBER = {379749},
MRREVIEWER = {H.\ Freudenthal},
}

@book {Kuczma,
    AUTHOR = {Kuczma, Marek},
     TITLE = {An introduction to the theory of functional equations and
              inequalities},
   EDITION = {Second},
      NOTE = {Cauchy's equation and Jensen's inequality,
              Edited and with a preface by Attila Gil\'{a}nyi},
 PUBLISHER = {Birkh\"{a}user Verlag, Basel},
      YEAR = {2009},
     PAGES = {xiv+595},
      ISBN = {978-3-7643-8748-8},
   MRCLASS = {39-02 (39Bxx)},
  MRNUMBER = {2467621},
       DOI = {10.1007/978-3-7643-8749-5},
       URL = {https://doi.org/10.1007/978-3-7643-8749-5},
}

@article {Console_Sergio_Fino_Evangelia,
    AUTHOR = {Console, Sergio and Fino, Anna and Samiou, Evangelia},
     TITLE = {The moduli space of six-dimensional two-step nilpotent {L}ie
              algebras},
   JOURNAL = {Ann. Global Anal. Geom.},
  FJOURNAL = {Annals of Global Analysis and Geometry},
    VOLUME = {27},
      YEAR = {2005},
    NUMBER = {1},
     PAGES = {17--32},
      ISSN = {0232-704X,1572-9060},
   MRCLASS = {22E25 (22E60 53C30)},
  MRNUMBER = {2130530},
MRREVIEWER = {Jir\B{o}\ Sekiguchi},
       DOI = {10.1007/s10455-005-2569-2},
       URL = {https://doi.org/10.1007/s10455-005-2569-2},
}

@incollection {Eberlein_moduli_2-step,
    AUTHOR = {Eberlein, Patrick},
     TITLE = {The moduli space of 2-step nilpotent {L}ie algebras of type
              {$(p,q)$}},
 BOOKTITLE = {Explorations in complex and {R}iemannian geometry},
    SERIES = {Contemp. Math.},
    VOLUME = {332},
     PAGES = {37--72},
 PUBLISHER = {Amer. Math. Soc., Providence, RI},
      YEAR = {2003},
      ISBN = {0-8218-3273-5},
   MRCLASS = {17B30},
  MRNUMBER = {2016090},
MRREVIEWER = {Vladimir\ V.\ Gorbatsevich},
       DOI = {10.1090/conm/332/05929},
       URL = {https://doi.org/10.1090/conm/332/05929},
}

@article {Homolya_Kowalski,
    AUTHOR = {Homolya, Szilvia and Kowalski, Old\v{r}ich},
     TITLE = {Simply connected two-step homogeneous nilmanifolds of
              dimension 5},
   JOURNAL = {Note Mat.},
  FJOURNAL = {Note di Matematica},
    VOLUME = {26},
      YEAR = {2006},
    NUMBER = {1},
     PAGES = {69--77},
      ISSN = {1123-2536,1590-0932},
      ISBN = {978-88-207-4018-4},
   MRCLASS = {53C30 (22E25)},
  MRNUMBER = {2267683},
MRREVIEWER = {Maura\ B.\ Mast},
}

@article {Jablonski_moduli,
    AUTHOR = {Jablonski, Michael},
     TITLE = {Moduli of {E}instein and non-{E}instein nilradicals},
   JOURNAL = {Geom. Dedicata},
  FJOURNAL = {Geometriae Dedicata},
    VOLUME = {152},
      YEAR = {2011},
     PAGES = {63--84},
      ISSN = {0046-5755,1572-9168},
   MRCLASS = {53C25 (17B30 22E25)},
  MRNUMBER = {2795236},
MRREVIEWER = {Megan\ M.\ Kerr},
       DOI = {10.1007/s10711-010-9546-z},
       URL = {https://doi.org/10.1007/s10711-010-9546-z},
}

@article {Dere,
    AUTHOR = {Der\'e, Jonas},
     TITLE = {The structure of underlying {L}ie algebras},
   JOURNAL = {Linear Algebra Appl.},
  FJOURNAL = {Linear Algebra and its Applications},
    VOLUME = {581},
      YEAR = {2019},
     PAGES = {471--495},
      ISSN = {0024-3795,1873-1856},
   MRCLASS = {17B30 (15A04 53C15)},
  MRNUMBER = {3989634},
MRREVIEWER = {Vesselin\ Drensky},
       DOI = {10.1016/j.laa.2019.07.032},
       URL = {https://doi.org/10.1016/j.laa.2019.07.032},
}

@book {Chevalley_Lie_Groups,
    AUTHOR = {Chevalley, Claude},
     TITLE = {Theory of {L}ie {G}roups. {I}},
    SERIES = {Princeton Mathematical Series},
    VOLUME = {vol. 8},
 PUBLISHER = {Princeton University Press, Princeton, NJ},
      YEAR = {1946},
     PAGES = {ix+217},
   MRCLASS = {20.0X},
  MRNUMBER = {15396},
MRREVIEWER = {G.\ Birkhoff},
}

@article {Levchuk,
    AUTHOR = {Levchuk, V. M.},
     TITLE = {Connections between the unitriangular group and certain rings.
              {II}. {G}roups of automorphisms},
   JOURNAL = {Siberian Mathematical Journal},
    VOLUME = {24},
      YEAR = {1983},
    NUMBER = {4},
     PAGES = {543--557},
}

@article {Graaf_char_not_2,
    AUTHOR = {de Graaf, Willem A.},
     TITLE = {Classification of 6-dimensional nilpotent {L}ie algebras over
              fields of characteristic not 2},
   JOURNAL = {J. Algebra},
  FJOURNAL = {Journal of Algebra},
    VOLUME = {309},
      YEAR = {2007},
    NUMBER = {2},
     PAGES = {640--653},
      ISSN = {0021-8693,1090-266X},
   MRCLASS = {17B30 (17B55)},
  MRNUMBER = {2303198},
MRREVIEWER = {Vladimir\ V.\ Gorbatsevich},
       DOI = {10.1016/j.jalgebra.2006.08.006},
       URL = {https://doi.org/10.1016/j.jalgebra.2006.08.006},
}

@article {Cicalo_Serena_Graaf,
    AUTHOR = {Cical\`o, Serena and de Graaf, Willem A. and Schneider, Csaba},
     TITLE = {Six-dimensional nilpotent {L}ie algebras},
   JOURNAL = {Linear Algebra Appl.},
  FJOURNAL = {Linear Algebra and its Applications},
    VOLUME = {436},
      YEAR = {2012},
    NUMBER = {1},
     PAGES = {163--189},
      ISSN = {0024-3795,1873-1856},
   MRCLASS = {17B30 (11E04 17B40 17B56)},
  MRNUMBER = {2859920},
MRREVIEWER = {Vesselin\ Drensky},
       DOI = {10.1016/j.laa.2011.06.037},
       URL = {https://doi.org/10.1016/j.laa.2011.06.037},
}

@article {Baer,
    AUTHOR = {Baer, Reinhold},
     TITLE = {Groups with abelian central quotient group},
   JOURNAL = {Trans. Amer. Math. Soc.},
  FJOURNAL = {Transactions of the American Mathematical Society},
    VOLUME = {44},
      YEAR = {1938},
    NUMBER = {3},
     PAGES = {357--386},
      ISSN = {0002-9947,1088-6850},
   MRCLASS = {20F99},
  MRNUMBER = {1501972},
       DOI = {10.2307/1989886},
       URL = {https://doi-org.ezproxy.lib.ou.edu/10.2307/1989886},
}

@misc{Demarche,
      title={A note on complex Lie Algebras isomorphic to their conjugate}, 
      author={Cyril Demarche},
      year={2026},
      eprint={2604.06979},
      howpublished={arXiv preprint, arXiv: 2604.06979},
      primaryClass={math.AG},
      url={https://arxiv.org/abs/2604.06979}, 
}

@misc{Borovoi-deGraaf-Guralnick,
      title={Constructing a complex Lie algebra isomorphic to its complex conjugate but not definable over reals}, 
      author={Mikhail Borovoi and Willem A. de Graaf and Robert M. Guralnick},
      year={2026},
      eprint={2607.19513},
      howpublished={arXiv preprint, arXiv:2607.19513},
      primaryClass={math.RA},
      url={https://arxiv.org/abs/2607.19513}, 
}

@article {Cornulier,
    AUTHOR = {Cornulier, Yves},
     TITLE = {Gradings on {L}ie algebras, systolic growth, and cohopfian
              properties of nilpotent groups},
   JOURNAL = {Bull. Soc. Math. France},
  FJOURNAL = {Bulletin de la Soci\'{e}t\'{e} Math\'{e}matique de France},
    VOLUME = {144},
      YEAR = {2016},
    NUMBER = {4},
     PAGES = {693--744},
      ISSN = {0037-9484,2102-622X},
   MRCLASS = {17B70 (20E07 20F18 20F40 20F69 20G15 22E25 22E40)},
  MRNUMBER = {3562610},
MRREVIEWER = {Laurent\ Bartholdi},
       DOI = {10.24033/bsmf.2725},
       URL = {https://doi-org.ezproxy.lib.ou.edu/10.24033/bsmf.2725},
}

\end{document}